\documentclass[11pt]{article}
\usepackage[letterpaper, margin=1.0in]{geometry}

\usepackage{microtype}
\usepackage{fullpage}
\usepackage[numbers,sort]{natbib}

\usepackage[utf8]{inputenc}							
\usepackage[T1]{fontenc}						

\usepackage{dsfont}							

\usepackage[
cal=cm,
]
{mathalfa}
\usepackage[labelfont={bf,small},labelsep=colon,font=small]{caption}	

\usepackage{hyperref}
\hypersetup{
final,
colorlinks=true,
linktocpage=true,
pdfstartview=FitH,
breaklinks=true,
pdfpagemode=UseNone,
pageanchor=true,
pdfpagemode=UseOutlines,
plainpages=false,
bookmarksnumbered,
bookmarksopen=false,
bookmarksopenlevel=1,
hypertexnames=true,
pdfhighlight=/O,
urlcolor=Maroon,linkcolor=MyBlue!60!black,citecolor=DarkGreen!70!black,	
pdftitle={},
pdfauthor={},
pdfsubject={},
pdfkeywords={},
pdfcreator={pdfLaTeX},
pdfproducer={LaTeX with hyperref}
}

\usepackage[dvipsnames,svgnames]{xcolor}						
\colorlet{MyBlue}{DodgerBlue!75!Black}
\colorlet{MyGreen}{DarkGreen!85!Black}

\renewcommand{\paragraph}{\subsection}

\usepackage{subcaption}

\usepackage{acronym}						
\usepackage{booktabs}						
\usepackage{paralist}						
\usepackage{xspace}							

\newcommand{\remove}[1]{{}}   
                 
\makeatletter
\def\paragraph{\@startsection{paragraph}{4}%
  \z@\z@{-\fontdimen2\font}%
  {\normalfont\bfseries}}
\makeatother

\usepackage[dvipsnames]{xcolor} 
\usepackage{tikz}
\usetikzlibrary{shapes}
\usetikzlibrary{positioning}
\usetikzlibrary{plotmarks}

\usepackage{microtype}
\usepackage{graphicx}
\usepackage{enumitem,caption}
\usepackage{booktabs} 
\usepackage{multirow}
\usepackage{adjustbox}

\usepackage{stmaryrd}
\usepackage{graphicx}

\usepackage{textcomp} 
\usepackage{quoting}
\usepackage{titletoc} 
\usepackage{fontawesome} 

\usepackage{algorithm}
\usepackage[noend]{algpseudocode}
\algrenewcommand\algorithmicrequire{\textbf{Input:}}
\algrenewcommand\algorithmicensure{\textbf{Output:}}

\renewcommand{\epsilon}{\varepsilon}

\newlength{\continueindent}
\usepackage{etoolbox}
\makeatletter
\newcommand*{\ALG@customparshape}{\parshape 2 \leftmargin \linewidth \dimexpr\ALG@tlm+\continueindent\relax \dimexpr\linewidth+\leftmargin-\ALG@tlm-\continueindent\relax}
\apptocmd{\ALG@beginblock}{\ALG@customparshape}{}{\errmessage{failed to patch}}
\makeatother

\algblock{ParFor}{EndParFor}
\algnewcommand\algorithmicpardo{\textbf{in parallel do}}
\algrenewtext{ParFor}[1]{\algorithmicfor\ #1\ \algorithmicpardo}
\algrenewtext{EndParFor}{\textbf{end for}}
\algtext*{EndParFor}
\algdef{SE}[SUBALG]{Indent}{EndIndent}{}{\algorithmicend\ }%
\algtext*{Indent}
\algtext*{EndIndent}

\usepackage{bm,mathtools,amsthm,amssymb}
\usepackage{bbm}

\usepackage{color}
\definecolor{puorange}{rgb}{0.80,0.20,0}
\definecolor{bluegray}{rgb}{0.04,0,0.7}
\definecolor{greengray}{rgb}{0.05,0.50,0.15}
\definecolor{darkbrown}{rgb}{0.40,0.2,0.05}
\definecolor{darkcyan}{rgb}{0,0.4,1}
\definecolor{black}{rgb}{0,0,0}
\definecolor{grey}{rgb}{0.93,0.93,0.93}

\newtheorem{theorem}{Theorem}
\newtheorem{lemma}[theorem]{Lemma}

\newtheorem{proposition}[theorem]{Proposition}
\newtheorem{corollary}[theorem]{Corollary}
\theoremstyle{definition}
\newtheorem{definition}[theorem]{Definition}

 \usepackage{thmtools}
 \declaretheoremstyle[
notefont=\bfseries, notebraces={}{},
bodyfont=\normalfont\itshape,
headformat=\NAME \NOTE
]{nopar}

\newcommand{\dl}[1]{\mathbb{#1}}
\newcommand{\fc}[2]{: #1 \rightarrow #2}
\newcommand{\prb}{\dl{P}}
\newcommand{\col}[1]{\left\{#1\right\}}
\newcommand{\mtr}{\mathsf{T}}
\newcommand{\galpha}{G_{\alpha}}
\newcommand{\reciprocal}[1]{\mathit{#1}^{-1}}
\newcommand{\invcon}[1]{\text{concave}-{#1}^{-1}}
\newcommand{\invconG}{\text{concave-}{G}^{-1}}
\newcommand{\invconcavityG}{\text{concavity-}{G}^{-1}}
\newcommand{\invconcavity}[1]{\text{concavity}-{#1}^{-1}}
\newcommand{\al}{\alpha}
\newcommand{\la}{\lambda}

\newcommand{\mtxt}[1]{\;\mbox{#1}\;}
\renewcommand{\Re}{\dl{R}}

\newcommand{\maps}{\mapsto}
\newcommand{\esp}[1]{\mathop{\rm \mathbb{E}\left (\mathit{#1} \right )}}
\newcommand{\foral}{\; \forall}

\newcommand{\scr}[1]{\mathcal{#1}} 
\newcommand{\sph}{\dl{S}^{m-1}}
\newcommand{\frk}{\mathfrak}
\newcommand{\hull}{\mathop{\rm Co}}
\newcommand{\dom}[1]{\mathrm{\mathcal{D}om}(#1)} 
\newcommand{\norm}[1]{\left\Vert#1\right\Vert}
\renewcommand{\phi}{\varphi}
\newcommand{\RR}{\mathbb{R}}
\newcommand{\C}{\mathcal{C}}

\theoremstyle{remark}
\newtheorem{remark}{Remark}						
\newtheorem*{remark*}{Remark}						
\newtheorem{example}{Example}						
\newtheorem*{example*}{Example}						

\title{On the Convexity of Level-sets of Probability Functions}

\author{
Yassine Laguel$^{1}$ \hspace{2em}
Wim Van Ackooij$^{2}$ \hspace{2em}
J\'{e}r\^{o}me Malick$^{1}$  \hspace{2em}
\vspace{1.em}
Guilherme Matiussi Ramalho$^{3}$ \\

\small{$^{1}$Univ. Grenoble Alpes, CNRS, LJK, 38000 Grenoble, France}  \\
\small{$^{2}$EDF R\&D, Saclay, France} \\
\small{$^{3}$Federal University of Santa Catarina, LABPLAN, Santa Catarina, Brazil}}
\date{\vspace{-2em}}
\begin{document}
\maketitle



\begin{abstract}
In decision-making problems under uncertainty, probabilistic constraints are a valuable tool to express safety of decisions. They result from taking the probability measure of a given set of random inequalities depending on the decision vector. Even if the original set of inequalities is convex, this favourable property is not immediately transferred to the probabilistically constrained feasible set and may in particular depend on the chosen safety level. In this paper, we provide results guaranteeing the convexity of feasible sets to probabilistic constraints when the safety level is greater than a computable threshold. Our results extend all the existing ones and also cover the case where decision vectors belong to Banach spaces. The key idea in our approach is to reveal the level of underlying convexity in the nominal problem data (e.g., concavity of the probability function) by auxiliary transforming functions. We provide several examples illustrating our theoretical developments.
\end{abstract}

\section{Introduction}

\subsection{Probability constraints and eventual convexity}

We consider a probabilistic constraint built up from the following
ingredients: a map $g \fc{X \times \Re^m}{\Re^k}$, where $X$ is a
(reflexive) Banach space, a random vector $\xi \in \Re^m$ defined
on an appropriate probability space, and a user-defined safety
level $p \in [0,1]$. The probabilistic constraint then reads:
\begin{equation}\label{eq:prb}
\phi(x) := \prb[ g(x,\xi) \leq 0 ] ~~\geq~~ p\,,
\end{equation}
where $\phi \fc{X}{[0,1]}$ is the associated probability function.
The interpretation of \eqref{eq:prb} is simple: one requires the
decision $x$ to be such that the random inequality system
$g(x,\xi) \leq 0$ holds with probability at least~$p$. Such
constraints (also called chance-constraint) often appear in
decision-making problems under uncertainty;  for the theory and
applications of chance-constraints optimization, we refer to
\cite{Prekopa_2003,Dentcheva_2009,Henrion_2004,vanAckooij_Henrion_Moller_Zorgati_2011,
Gonzalez-Gradon_Heitsch_Henrion_2017} and references therein.

In this paper we focus on the ``convexity of the probabilistic
constraint" \eqref{eq:prb}, i.e., the convexity of the set of
feasible solutions defined by
\begin{equation}\label{eq:Mp}
M(p) := \col{ x \in X\; : \prb[ g(x,\xi) \leq 0 ] \geq p }.
\end{equation}
Understanding when $M(p)$ is a convex set is important for the
point of view of optimization, to guarantee that local solutions
are also globally optimal and to use numerical solution methods
that exploit this convexity (we review the most popular methods in
section~\ref{sec:methods}). A first result of the convexity of
$M(p)$ follows from Pr{\'e}kopa's celebrated log-concavity theorem
(see \cite[Proposition 4]{Farshbaf-Shaker_Henrion_Homberg_2017}
for its infinite dimensional version and \cite{Diniz_Henrion_2017}
for generalizations):
the convexity of $M(p)$ is guaranteed  for all~$p\in [0,1]$, when
$-g$ is jointly quasi-concave in both arguments and $\xi$ an
appropriate random vector. However, joint-quasi-concavity of $-g$
is rather ``exceptional" and fails in many basic situations. For
example, when $g(x,\xi) = x^\mtr \xi$ and $\xi$ is multi-variate
Gaussian, it is well known that $M(p)$ is convex only whenever $p
\geq \frac 12$ (see e.g., \cite{Kataoka_1963}). In this example
and many others, we thus observe that if the convexity of $M(p)$
does not hold for all $p$, there still exists a (computable)
threshold $p^\ast \in [0,1]$ such that the set $M(p)$ is convex
for all $p \geq p^\ast$. This property is called \emph{eventual
convexity}
as observed by \cite{Prekopa_2001} and coined by \cite{Henrion_Strugarek_2008} (which studies the
case where $g$ is separable and $\xi$ has independent components).
Eventual convexity results are further generalized in
\cite{Henrion_Strugarek_2011} by allowing for the components of
$\xi$ to be coupled through a copul\ae\/ dependency structure.
These results are refined, by allowing for more copul\ae\/ and
with sharper bounds for $p^\ast$ in \cite{vanAckooij_2013}, and
extended to all Archimedian copul\ae\/ in
\cite{vanAckooij_Oliveira_2016}, where also an appropriate
solution algorithm is provided. When the mapping $g$ is
non-separable, eventual convexity results are provided in
\cite{vanAckooij_Malick_2017} for the special case where $\xi$ is
elliptically symmetrically distributed. Here we will simplify,
clarify and extend these results.

\subsection{Ideas, contributions, and outline of this paper}

In this paper, we build on this line of research about
establishing convexity of superlevel-sets of probability functions
\eqref{eq:Mp} for $p$ larger than a threshold. We show that a
notion of
generalized concavity 
naturally appears in this framework, and allows us to reveal the
level of hidden convexity of the data. We formalize a way to
analyze separately the convexity inherent to the randomness and
the one associated with the optimization model structure.

Roughly speaking, our approach is the following. In various
contexts, the probability function involves a composition of two
functions $F \circ Q$, with $F\colon\RR\rightarrow[0,1]$ carrying
the randomness of the problem and $Q\colon\RR^n\rightarrow\RR$
given by the optimization model. We split the problem of
establishing concavity (or at least quasi-concavity) of the
composition $F \circ Q$  by finding an adequate function $G$ (the
inverse of which is denoted $\reciprocal{G}$) to write
\begin{equation}\label{eq:compo}
F \circ Q = F \circ \reciprocal{G} \circ G \circ Q
\end{equation}
such that  $F \circ \reciprocal{G}$ and $G \circ Q$ satisfy
appropriate convexity properties. Thus this approach naturally
raises interest in the concavity of function $G \circ Q$, which is
a ``transformable concavity", formalized by the notion called
transconcavity or $G$-concavity by \cite{Tamm_1977} and
\cite{Avriel_Diewert_Schaible_Zang_2010}. Similarly, we will
briefly study concavity of terms $F \circ \reciprocal{G}$ and
introduce the counterpart notion that we will call
$\invconcavityG$.

Beyond this intuition, the analysis of the convexity of
chance-constrained sets is not trivial: $F$, as a distribution
function, cannot be concave and one has to be careful with working
on appropriate smaller subsets. The size of those subsets turns
out to be directly related to the level of probability beyond
which convexity is guaranteed. Arguments of this type are
implicitly and partially used in \cite{Henrion_Strugarek_2008,
vanAckooij_2013,vanAckooij_Oliveira_2016,vanAckooij_Malick_2017}.
Here we highlight this approach, exploit it in its full generality, and provide a set of tools to apply it in practice. We apply this set of tools to get eventual convexity statements in two general contexts
\begin{itemize}
\item[(A)] when $\xi$ is specific (elliptically distributed) and
$g$ is general (with light geometrical assumptions),
\smallskip
\item[(B)] when $\xi$ is general (with known associated copula)
and $g$ is specific (of the form $g(x,\xi)=\xi-h(x)$)
\end{itemize}
In these two cases, we prove the existence of the threshold $p^*$
such that $M(p)$ is convex for all $p\geq p^*$. Beyond the
theoretical contribution regarding the existence of $p^*$, we also
provide a concrete way of numerically evaluating $p^*$ from the
nominal data. Our results are illustrated through various examples
(not captured by existing results) of eventual convexity with a
specified threshold $p^*$.

The contributions of this work are thus the following ones. Our
main contribution is to identify the interplay of generalized
concavity with the functions used to model the constraints and the
uncertainty in chance-constrained optimization problems. This
clarification would allow practitioners to refine their models of
nonlinearity and readily swapping certain uncertainty
distributions by others sharing similar generalized concavity
properties. Since we care about practical use of our results, this
introduction is meant to be an accessible overview of the
state-of-the-art on eventual convexity and its use in practice. We
also provide througout the paper many examples illustrating our
results. Finally we specifically list our main technical
contributions:
\begin{itemize}
\item We slightly extend the existing notion of $G$-concavity (by
allowing for decreasing functions $G$) and introduce the new
notion of $\invconcavityG$ that goes with it. We also extend a
useful result of \cite{Henrion_Strugarek_2008} on one-dimensional
distribution functions and give necessary and sufficient
conditions under which such functions can be composed with
monotonous maps to make them concave.

\medskip

\item For the context (A), our results clarify and extend those of
\cite{vanAckooij_Malick_2017} allowing us to treat new situations:
we use various forms of $G$-concavity and thus can cover a wider
range of non-linear mappings $g$.

\medskip

\item For the context (B), we first refine the results of
\cite{Zadeh_Khorram_2012} by providing a better threshold. More
importantly, our results cover new situations as they tackle the
most general case with nonlinear mappings and copul\ae, while
\cite{vanAckooij_2013,vanAckooij_Oliveira_2016} restricts to
$\alpha$-concavity (a special $G$) and \cite{Zadeh_Khorram_2012}
restricts to independent copul\ae.

\medskip

\item We also extend all these previous results (that consider
finite dimensional decision vectors $x$), by analysing and stating
our results in Banach spaces. This opens the door to cover recent
applications in PDE-constrained optimization (e.g.,
\cite{Farshbaf-Shaker_Henrion_Homberg_2017}).
\end{itemize}

%

The paper is organized as follows. After section
\ref{sec:methods}, which will put this work into a broader
practical perspective, our development starts in Section
\ref{sec:genconc}, which provides a careful account of useful and
used notions of generalized concavity together with calculus
rules. Section \ref{sec:interplay} then provides general eventual
convexity statements in the two contexts (A) and (B) in
Sections\;\ref{subsec:nonlin} and\;\ref{subsec:copulae}
respectively.
Finally Section \ref{sec:examples} is devoted to providing
examples covered in our extended framework.

\subsection{Discussion on applicability: practicality \& open issues}
\label{sec:methods}

Although this work is theoretical, it can be inserted in the
bigger picture of solving chance-constrained optimization
problems. Of course, a priori, evaluating probability functions is
computationally demanding (see e.g., the discussion in
\cite{Nemirovski_Shapiro_2006a}), especially if the random vector
is highly dimensional. Still recent results (e.g.,
\cite{Bremer_Henrion_Moller_2015}) indicate that random vectors
with dimensions in their hundreds, i.e., practically relevant
sizes, can be handled with CPU times hovering roughly around a
minute. Those CPU times relate to solving a non-convex
optimization problem involving a probability function, evaluated
repeatedly. Although this does not alleviate the theoretical or
algorithmic difficulties, it does show that, by exploiting the
structure present in applications, one can solve applications of
relevant size. The existing numerical solution methods use
sample-based approximations of the probabilistic constraint or
treat probability functions (or a surrogate) by non-linear
optimization techniques; we briefly review the main methods here.
Notice that many of these methods rely on some convexity
properties of the chance-constrained problems, bringing interest
to the results of this paper.

Popular numerical methods for dealing with 
probabilistic constraints are sample-based approximations, e.g.,
\cite{Luedtke_Ahmed_2008,Pagnoncelli_Ahmed_Shapiro_2009,Luedtke_2010,Luedtke_Ahmed_Nemhauser_2010}
with various strenghtening procedures, e.g.,
\cite{Kucukyavuz_2012,Wu_Kucukyavuz_2019} or investigation of
convexification procedures \cite{Ahmed_Xie_2018}. We can also
mention boolean approaches, e.g.,
\cite{Lejeune_2012,Lejeune_Margot_2016,Kogan_Lejeune_2014,Kogan_Lejeune_Luedtke_2016},
$p$-efficient point based concepts, e.g.,
\cite{Dentcheva_Prekopa_Ruszczynski_2000,Lejeune_Noyan_2010,Dentcheva_Martinez_2012b,Dentcheva_Martinez_2013,vanAckooij_Berge_Oliveira_Sagastizabal_2017},
robust optimization \cite{Ben-Tal_Nemirovski_2000},  penalty
approach \cite{Ermoliev_Ermolieva_Macdonald_Norkin_2000}, scenario
approximation \cite{Calafiore_Campi_2006,Ramponi_2018}, convex
approximation \cite{Nemirovski_Shapiro_2006a}, or yet other
approximations
\cite{Hong_Yang_Zhang_2011,Geletu_Hoffmann_Kloppel_Li_2015}. Aside
from this rich literature, the non-linear constraint
\eqref{eq:prb} can also be dealt with directly as such, from the
study of (generalized) differentiability of probability functions
and the development of readily implementable formul\ae\/ for
gradients. Such formul\ae\/ can be further improved by using well
known ``variance reduction'' techniques, such as Quasi-Monte Carlo
methods (e.g., \cite{Brauchart_Saff_Sloan_Womersley_2014}) or
importance sampling (e.g.,
\cite{Barrera_Homem-de-Mello_Moreno_Pagnoncelli_Canessa_2016}).
For further insights on differentiability, we refer to
e.g.,\,\cite{vanAckooij_Henrion_2016,Henrion_Moller_2012,Uryasev_1995,Uryasev_2009,Kibzun_Uryasev_1998,Royset_Polak_2007,Marti_1995}.
Nonlinear programming methods using these properties include
sequential quadratic programming \cite{Bremer_Henrion_Moller_2015}
and the promising bundle methods
\cite{vanAckooij_Sagastizabal_2014,vanAckooij_Oliveira_2016}.

Practical probabilistic constrained problems also involve several
other constraints, that can be represented as an abstract subset
$S\subset X$. Important questions concern, in fact, the
constrained set $M(p)\cap S$, for which the results presented in
this paper might be used. As a brief observation we do write $p$
and not $p^*$, since $p$ is the user chosen safety level and thus
what is practically relevant.
\begin{itemize}
\item Is $M(p)\cap S$ convex? Convexity of $S$ is achieved in many
practical cases: in a significant share of applications $S$ is
polyhedral, or easily seen to be convex. The difficulty in
establishing convexity of $M(p)\cap S$ therefore lies in checking
whether $M(p)$ is a convex set, which is the aim of this paper. As
a sufficient condition, the user would check $p \geq p^*$, when
applicable, in order to have a global guarantee on computed
solutions. Since, occasionally, $p^*$ may depend adversely on
random vector dimension, if the test fails, this does not
necessarily imply that $M(p) \cap S$ is not convex. The results of
this paper rather indicate that a computed ``solution'' should not
necessarily be taken as a global solution, because convexity of
$M(p) \cap S$ is no longer guaranteed. Such an information is
still useful for the user, who may decide to invest additional
effort in running the local optimization solver with multiple
starting points, or calling another (expensive) global
optimization solver.
\smallskip
\item Is $M(p)\cap S$ non-empty ? The safety level $p$ is chosen
by the user, who, as a modeler, is responsible for ensuring that a
well-posed model is formulated. In practice, we can expect a
reasonable convex optimization solver to return an infeasibility
flag when the set $M(p)\cap S$ is empty (in the case that $M(p)
\cap S$ is convex as attested by the previous point). In such a
case, the user can examine his data, and subsequently formulate a
better model. Answering the feasibility question without convexity
is of course an entirely different matter, and such a
theoretically and algorithmically difficult problem goes largely
beyond the scope and the setting of this paper. Below we mention
some relevant heuristic procedures that have worked well in our
experience.
\end{itemize}
The feasibility regarding probabilistic restrictions is
related to the question of maximizing the probability function
over $X$ or $S$, which recently has received special attention;
see e.g.,
\cite{Fabian_Csizmas_Drenyovszki_vanAckooij_Vajnai_Kovacs_Szantai_2018,Fabian_Csizmas_Drenyovszki_Vajnai_Kovacs_Szantai_2019,Minoux_Zorgati_2016}.
Indeed, if the maximal probability thus found is greater than or
equal to $p$, then $M(p)\cap S$ is ensured to be a nonempty set.
Of course, finding the global solution of this probability
maximization problem is a hard task in general, because the
probability function need not be concave. 
Heuristically, feasibility can also be addressed by considering a
sample based variant of probabilistic maximization problem, with
few samples, i.e.,
\begin{align*}
\min_{z_1,...,z_N}\quad & \sum_{i=1}^N z_i \\
\mbox{s.t.}\quad & g(x,\xi_i) \leq Mz_i, \quad x \in S, ~z_i \in \col{0,1},
\end{align*}
where $M$ is an appropriate ``big-M'' constant, $S$ is the
deterministic constraints set, and $\xi_1,...,\xi_N$ are i.i.d.
samples of\;$\xi$. The last problem can in principle be solved with
fairly few samples (small $N$) and with low accuracy (e.g., 10\%
MIP-gap) to produce $\bar x$. An a posteriori evaluation of the
probability function can then assert feasibility of $\bar x$ for
the true probabilistic constraint. Indeed, when $g$ is convex in
$x$ and $S$ is convex, one can solve the last program with the
methodology laid out in \cite{vanAckooij_Frangioni_Oliveira_2016}.
The sample can also be exploited in a ``scenario approach'' (the asymptotics with respect to $N$ are well studied in
for instance \cite{Ramponi_2018}). Finally, the approach
(maximization of copula structured probability not requiring
convexity) in \cite{vanAckooij_Oliveira_2019} can also be
employed. It consists of minimizing a lower-$C^2$ function
(requiring easily verified differentiability, whenever the copula
is Archimedian) with tools from nonsmooth optimization.


\section{Generalized concavity, propagation of concavity, and cumulative distribution functions}
\label{sec:genconc}


In this section, we gather the tools on generalized concavity that
we will use in next sections to reveal the underlying concavity of
nominal data in probabilistic constraints. Section
\ref{sec:transconc} briefly reviews the definitions and useful
properties of $G$-concavity (also called transconcavity (see
\cite{Avriel_Diewert_Schaible_Zang_2010})) and Section
\ref{sec:concaveg} introduces the right counterpart of
$\invconcavityG$. We provide new technical lemmas, including a
characterization of generalized concavity of cumulative
distribution functions.

\subsection{Discussions on \texorpdfstring{$G$-}concavity}\label{sec:transconc}

We start by recalling the notion of $G$-concavity introduced by
\cite{Tamm_1977} and presented in the book
\cite{Avriel_Diewert_Schaible_Zang_2010} under the name
transconcavity. We just add here the possibility of $G$ being
strictly decreasing, which will turn out to be useful in our
context.

\begin{definition}
Let $X$ be a Banach space and $C$ be a convex subset of $X$. We
say that a function $f\colon C \rightarrow \mathbb{R}$ is
$G$-concave if there exists a continuous and strictly
monotonic\footnote{In this paper, we denote the inverse of a map $G$
by $\reciprocal{G}$, and the division by a map, whenever well
defined, by $\frac{1}{G}$.} function $G : f(C) \rightarrow
\mathbb{R}$ such that
\begin{equation*}
f(\lambda x + (1-\lambda) y) \geq \reciprocal{G}(\lambda G\circ
f(x) + (1-\lambda) G\circ f(y))
\end{equation*}
holds for all $x,y \in C$ and $\lambda \in [0,1]$.
\label{def:transconcave}
\end{definition}

Note that when $G$ is increasing, $G$-concavity of $f$ is just the
concavity of the map $G \circ f$. When $G$ is decreasing,
$G$-concavity of $f$ is simply the convexity of $G \circ f$. A
given function $f$ can be ``$G$-concave'' for several different
mappings $G$. It will be convenient, however, to pin down a specific
choice and subsequently speak of $G$-concavity of $f$ for such a
specific choice. The naming ``transconcavity" would then refer to
an unspecified, yet implicitly assumed to exist, mapping $G$ for
which $f$ is $G$-concave.

\begin{example}[Special family]
A particularly well studied set of choices for $G$ is that of the
family
\begin{equation}\label{eq:galpha}
\galpha : t \mapsto t^\alpha \quad \text{for $\alpha \in
\Re\setminus\col{0}$} \qquad \text{and} \qquad G_0 : t \mapsto
\ln(t).
\end{equation}
This family has several properties that help to measure a ``level
of generalized concavity'' of a function $f$, as used below in the
definition of $\al$-concavity.\label{ex:galpha} \hfill\qedsymbol
\end{example}


We introduce the following mapping $m_{\alpha} \fc{\dl{R}_+ \times
\dl{R}_+ \times [0,1]}{\dl{R}}$ (for a given $\alpha \in
[-\infty,\infty)$) defined as follows:
\begin{equation}
\mtxt{if} ab = 0 \mtxt{and} \alpha \leq 0, \qquad
m_{\alpha}(a,b,\lambda) = 0 \label{eq:ma1}
\end{equation}
else, for $\lambda \in [0,1]$, we let:
\begin{equation}
m_{\alpha}(a,b,\lambda) = \left \{
\begin{array}{ccc}
a^{\lambda}b^{1-\lambda} &\mtxt{if}& \alpha=0 \\
\min\col{a,b} &\mtxt{if}& \alpha=-\infty \\
(\lambda a^{\alpha} + (1-\lambda)b^\alpha)^{\frac{1}{\alpha}}
&\mtxt{else}&
\end{array}
\right.
\end{equation}

This enables us to define the known notion of $\alpha$-concavity
as a ``particular case'' of $G$-concavity.
\begin{definition}[$\alpha$-concave function] \label{def:convfunction}
Let $X$ be a Banach space and $C$ be a convex subset of $X$. We
say that a function $f : C \rightarrow \mathbb{R}_+$ is
$\alpha$-concave if
\begin{equation}
f(\lambda x + (1-\lambda) y) \geq m_{\alpha}(f(x),f(y),\lambda),
\end{equation}
for all $x,y \in C$ and $\lambda \in [0,1]$.
\end{definition}

Observe that for $\alpha \neq 0$, $\alpha$-concavity of $f$ is
indeed equivalent with $\galpha$-concavity of $f$ (in the sense of
Definition\;\ref{def:transconcave} with $\galpha$ of
\eqref{eq:galpha}). Notice that $1$-concavity coincides with the
usual notion of concavity.

We also note that our definition of the mapping $m_{\alpha}$
differs slightly of that found in \cite[Def 4.7]{Dentcheva_2009},
in so much that we have appended the condition $\alpha \leq 0$ to
condition \eqref{eq:ma1}. The reason for this is that otherwise
the definition does not match up with what is expected whenever
$a=0$ or $b=0$ and $\alpha > 0$. In particular consider $\alpha=1$
and the usual definition of concavity for a function $f$, two
points $x,y \in C$ with $f(y)=0$ for instance. Since our
definition slightly differs from the classical one, we provide the
proof of the following technical lemma used to establish the
hierarchy of $\alpha$-concavity.


\begin{lemma}
Let $a, b \in \Re_+$, $\lambda \in [0,1]$ be given and fixed. Then
for $\alpha, \beta \in [-\infty,\infty)$, $m_\beta(a,b,\lambda)
\leq m_\alpha(a,b,\lambda)$ holds  when $\beta \leq \alpha$.
Moreover the map $\alpha \maps m_\alpha(a,b,\lambda)$ is
continuous. \label{lem:hierarchy}
\end{lemma}

\begin{proof}
Let $a,b,\lambda$ be as in the statement. Since $\lambda$ is
arbitrary and $m_{\alpha}(a,b,\lambda) = m_{\alpha}(b,a,1 -
\lambda)$, we may without loss of generality assume $a \leq b$.
Furthermore, let also $\beta \leq \alpha$ be given but fixed. We
will proceed by a case distinction.

\noindent \textit{First case}: If $\alpha \geq \beta > 0$ or
$\alpha
> 0 > \beta$ then the mapping $t \mapsto t^{\frac{\alpha}{\beta}}$
is convex on $\RR_+$. So, we have
\begin{equation*}
\lambda u^{\frac{\alpha}{\beta}} +
(1-\lambda)v^{\frac{\alpha}{\beta}} \geq (\lambda u +
(1-\lambda)v)^{\frac{\alpha}{\beta}}, \qquad \text{for $u,v\geq
0$}
\end{equation*}
which since $\alpha > 0$, and hence $t \maps t^{\frac{1}{\alpha}}$
strictly increasing on $\RR_+$, is equivalent to $(\lambda
u^{\frac{\alpha}{\beta}} +
(1-\lambda)v^{\frac{\alpha}{\beta}})^{\frac{1}{\alpha}} \geq
(\lambda u + (1-\lambda)v)^{\frac{1}{\beta}}$. The desired result
follows by substituting $u = a^\beta$ and $v = b^\beta$.

\noindent \textit{Second case}: If $0 > \alpha \geq \beta$, then
we have $-\beta \geq - \alpha > 0$ and for any $u,v\geq 0$, we can
apply the previous case to obtain the inequality:
\begin{equation*}
(\lambda u^{- \beta} + (1-\lambda) v^{-\beta})^{\frac{-1}{\beta}}
\geq (\lambda u^{-\alpha} + (1-\lambda)
v^{-\alpha})^{\frac{-1}{\alpha}}.
\end{equation*}
The latter is equivalent with
$$(\lambda (\frac{1}{u})^{\beta} + (1-\lambda) (\frac{1}{v})^{\beta})^{\frac{1}{\beta}} \leq (\lambda (\frac{1}{u})^{\alpha} + (1-\lambda)
(\frac{1}{v})^{\alpha})^{\frac{1}{\alpha}},$$ provided that $u,v >
0$ hold. Assuming $a > 0$ (and hence $b>0$), we may substitute $u
= \frac{1}{a}$ and $v = \frac{1}{b}$ to obtain the desired
inequality. When $ab=0$, since both $\alpha$ and $\beta \leq 0$,
the desired inequality holds.
%

\noindent\textit{Third case}: To treat a case where $\alpha = 0$
or $\beta=0$, we first establish continuity of $\alpha \mapsto
m_\alpha$ around 0. To this end, consider the following Taylor
expansions:
\begin{align*}
a^\alpha &= e^{\alpha \ln(a)} = 1 + \alpha \ln(a) + o(\alpha) \\
\frac{1}{\alpha} \ln(\lambda a^\alpha + (1-\lambda) b^\alpha) &=
\frac{1}{\alpha} \ln(1 + \alpha(\lambda \ln(a) + (1-\lambda)
\ln(b)) + o(\alpha))
= \lambda \ln(a) + (1-\lambda) \ln(b) + o(1)
\end{align*}
Consequently, we get at the limit when $\alpha \rightarrow 0$
\begin{equation*}
\exp(\frac{1}{\alpha} \ln(\lambda a^\alpha + (1-\lambda)
b^\alpha)) = \exp(\lambda \ln(a) + (1-\lambda) \ln(b) + o(1)))
\rightarrow a^\lambda b^{1-\lambda} = m_0(a,b,\lambda),
\end{equation*}
Now assume $\beta=0$ and $\alpha
> 0$. If $m_{0}(a,b,\lambda) \leq m_{\alpha}(a,b,\lambda)$ were
not to hold, then it would follow that
$m_{\alpha}(a,b,\lambda) < m_{0}(a,b,\lambda)$. We may pick a
sequence $\alpha_k \downarrow 0$ and for $k$ large enough it holds
$0 < \alpha_k \leq \alpha$. By the already established order, we
have $m_{\alpha_k}(a,b,\lambda) \leq m_{\alpha}(a,b,\lambda) <
m_{0}(a,b,\lambda)$. But this contradicts the just established
continuity of $\alpha \mapsto m_{\alpha}(a,b,\lambda)$ near $0$.
The situation $\alpha=0$ can be established along similar lines of
argument.
%

\noindent \textit{Fourth case}: The situation wherein $\beta < 0 <
\alpha$ follows by invoking the third case twice, since indeed
$m_{\beta}(a,b,\lambda) \leq m_{0}(a,b,\lambda) \leq
m_{\alpha}(a,b,\lambda)$.

\noindent
\textit{Last case}: When $\beta=-\infty$. The inequality
$m_{-\infty}(a,b,\lambda) \leq m_{\alpha}(a,b,\lambda)$ holds
trivially whenever $ab=0$. By combining the previous cases we may
assume $\alpha < 0$ as well as $a, b > 0$. We observe that $a,b
\geq m_{-\infty}(a,b,\lambda)=\min\col{a,b}$ and since $\alpha <
0$, it holds $a^\alpha, b^\alpha \leq
m_{-\infty}(a,b,\lambda)^{\alpha}$. Consequently too,
$$\lambda a^\alpha + (1-\lambda) b^\alpha \leq
m_{-\infty}(a,b,\lambda)^{\alpha}.$$ Since $t \maps
t^{\frac{1}{\alpha}}$ is strictly decreasing, we get
$m_{-\infty}(a,b,\lambda) \leq m_{\alpha}(a,b,\lambda)$.

We finish by proving the continuity. Since both terms $\lambda
a^{\alpha}$ and $(1-\lambda) b^{\alpha}$ are nonnegative, we have
that:
\begin{equation*}
\max(\lambda a^{\alpha}, (1-\lambda) b^{\alpha}) \leq \lambda
a^{\alpha} + (1-\lambda) b^{\alpha},
\end{equation*}
 and consequently
\begin{equation}\label{eq:sandwhich}
m_{-\infty}(a,b,\lambda) \leq (\lambda a^{\alpha} + (1-\lambda)
b^{\alpha})^{\frac{1}{\alpha}} \leq \max(\lambda a^{\alpha},
(1-\lambda) b^{\alpha})^{\frac{1}{\alpha}}
=[\max(\lambda^{\frac{1}{-\alpha}} \frac{1}{a},
(1-\lambda)^{\frac{1}{-\alpha}} \frac{1}{b})]^{-1}.
\end{equation}
Hence by passing to the limit:
\begin{equation*}
\lim_{\alpha \rightarrow -\infty}
[\max(\lambda^{\frac{1}{-\alpha}} \frac{1}{a},
(1-\lambda)^{\frac{1}{-\alpha}} \frac{1}{b})]^{-1} =
[\max(\frac{1}{a}, \frac{1}{b})]^{-1} = \min(a,b),
\end{equation*}
This gives the continuity at $-\infty$.
\end{proof}

The above property of the map $m_{\alpha}$ allows us to establish
an entire hierarchy of ``concavity'' immediately, as formalized by
the next corollary.

\begin{corollary}[Hierarchy of $\alpha$-concavity]
Let $X$ be a Banach space and $C$ be a convex subset of $X$. Let
the map $f : C \rightarrow \mathbb{R}_+$, together with
$\alpha,\beta \in [-\infty,\infty)$, be given. If $f$ is
$\alpha$-concave, it is also $\beta$-concave when $\alpha \geq
\beta$. In particular $f$ is quasi-concave. \label{cor:alpconc}
\end{corollary}

The family of mappings $\col{\galpha}_{\alpha}$ of
\eqref{eq:galpha} allows us to distinguish the level of
generalized concavity of a function $f$ ranging from
quasi-concavity ($\alpha = -\infty$) to classic concavity ($\alpha
= 1$). Intuitively, the greater\;$\alpha$ is, the ``more
concave'', $f$ will be.
All the practical examples in this paper will use these functions
to quantify and extract underlying convexity. Let us mention
though that this family does not capture completely the subtle
notion of transconcavity (see Example\;\ref{ex:tbigalp} below) and
that alternative families of functions could be considered, such
as the exponential family of functions $G_r : t \mapsto -e^{-rt}$,
for varying values of $r$, extensively studied in
\cite[Chap.8]{Avriel_Diewert_Schaible_Zang_2010}.

\begin{example}[transconcavity does not imply $\alpha$-concavity]
\label{ex:tbigalp} Let us provide an example of a mapping $h
\fc{\Re}{\Re_+}$ that is not $\alpha$-concave for any $\alpha \in
\Re$, but is $G$-concave for an appropriate choice of a map $G$.
We will show that $h \fc{\Re}{\Re_+}$ defined as $h(x) =
\exp(-x^3)$ is not $\alpha$-concave for any $\alpha < 0$ (and then
by Corollary\;\ref{cor:alpconc} that $h$ can not be
$\alpha$-concave for any $\alpha$). Indeed, let $\alpha < 0$ be
arbitrary. Then $\alpha$-concavity of $h$ is equivalent to
convexity of $h^{\alpha}$. Now by differentiating twice we obtain:
\begin{align*}
\frac{d}{dx} (h^\alpha)(x) &= -3x^2 \alpha (h^\alpha)(x)
\qquad\text{and}\qquad \frac{d^2}{dx^2} (h^\alpha)(x) = h^{\alpha}(x)
\alpha x (9 \alpha x^3 - 6 ).
\end{align*}
For $x < 0$, $\alpha x > 0$ and moreover $x \maps 9 \alpha x^3 -
6$ has unique (negative root) $x^r = (\frac{6}{9\alpha})^{\frac
13}$. Consequently $9 \alpha x^3 - 6 > 0$ for $x < x^r$ and $9
\alpha x^3 - 6 < 0$ for $x > x^r$. By combining with the above, we
establish that $\frac{d^2}{dx^2} (h^\alpha)(x) < 0$ must hold for
$x \in (x^r, 0)$, implying that $h^\alpha$ can not be convex,
i.e., $h$ is not $\alpha$-concave.

Now define $G \fc{(0,+\infty)}{\Re_+}$ as $G(x) = \exp(-(\ln(x))^{\frac
13})$, then by direct computation $\frac{d}{dx} G(x) = -
\frac{1}{3x (\ln(x))^{\frac 23} }G(x)$ at any $x \neq 1$. We now
readily verify that $G$ is strictly decreasing. Moreover, $G(h(x))
= \exp(x)$; hence, by definition, $h$ is $G$-concave.\hfill\qedsymbol
\end{example}

We now recall and extend a lemma from
\cite{Avriel_Diewert_Schaible_Zang_2010} that enables us to
propagate the property of $G$-concavity.

\begin{lemma}[Propagation of generalized concavity]
Let $X$ be a Banach space and $C$ be a convex subset of $X$. Let
the map $f \fc{X}{\Re}$ be a $G_1$-concave function for an
appropriate choice $G_1 \fc{f(C)}{\Re}$. Let $G_2$ be a continuous
and strictly monotonic function over $f(C)$. If
$\reciprocal{G_1}$, the inverse function of $G_1$ is $G_2$-concave
over $G_1\circ f(C)$, then $f$ is also $G_2$-concave over $C$.
\label{lem:G1impliesG2}
\end{lemma}

\begin{proof}
By assumption we have, for any $x,y \in C, \lambda \in [0,1]$ that
\begin{equation*}
f(\lambda x + (1-\lambda) y) \geq \reciprocal{G_1}(\lambda
G_1\circ f(x) + (1-\lambda) G_1\circ f(y))
\end{equation*}
and for any $u,v \in G_1\circ f(C), \lambda \in [0,1]$, we have
\begin{equation*}
    \reciprocal{G_1}(\lambda u + (1-\lambda) v) \geq \reciprocal{G_2}(\lambda G_2\circ \reciprocal{G_1} (u) + (1-\lambda) G_2\circ
    \reciprocal{G_1}(v)).
\end{equation*}

Hence, if we fix $x,y \in C, \lambda \in [0,1]$ and we set $u =
G_1\circ f(x), v = G_2\circ f(y)$, we get from these two
inequalities:
\begin{align*}
f(\lambda x + (1-\lambda) y) &\geq \reciprocal{G_1} (\lambda u +
(1-\lambda) v)
                             \geq \reciprocal{G_2}(\lambda G_2\circ \reciprocal{G_1}(u) + (1-\lambda) G_2\circ \reciprocal{G_1}(v))\\
                             & = \reciprocal{G_2}(\lambda G_2\circ \reciprocal{G_1}(G_1\circ f(x)) + (1-\lambda) G_2\circ \reciprocal{G_1}(G_1\circ f(y))) \\
                             & \geq \reciprocal{G_2}(\lambda G_2\circ f(x) + (1-\lambda) G_2\circ f(y)),
\end{align*}
which gives the result.
\end{proof}

\subsection{Study of concavity-\texorpdfstring{$G^{-1}$}.}\label{sec:concaveg}


We introduce in this section the notion of $\invconcavityG$ which
is the right counterpart of the classical $G$-concavity recalled
previously. In view of \eqref{eq:compo}, the two complementary
notions will be useful in the sequel, in particular because we
establish that many cumulative distribution functions are
$\invconG$. Along the way, we generalize a result of
\cite{Henrion_Strugarek_2008}.

\begin{definition}[$\invconG$ functions]\label{def:concave-g}
Let $F \fc{\Re}{\Re}$ and $G \fc{\Re}{\Re}$ be continuous and
strictly monotonic mappings. The map $F$ is said to be $\invconG$
on an interval $I \subset \Re$ if $F \circ \reciprocal{G}$ is
concave on the interval $I$. By extension of Definition
\ref{def:concave-g}, we also speak of a $G_1$-concave-$G_2^{-1}$
function $F$ if $G_1 \circ F \circ \reciprocal{G_2}$ is concave.
\end{definition}

This definition can be specialized as follows when considering the
family $\col{\galpha}_{\alpha}$ of Example\;\ref{ex:galpha}. Let
$\alpha \in (-\infty, 1]$ be given; we say that $f:\mathbb{R}
\rightarrow \mathbb{R}$ is concave-$\alpha$ on an interval $I
\subset \Re$ if (i) $f$ is increasing and $t \mapsto
f(t^{\frac{1}{\alpha}})$ is concave on $I$, or (ii) $f$ is
decreasing and $t \mapsto f(t^{\frac{1}{\alpha}})$ is convex on
$I$. (Note that $t \mapsto f(\exp(t))$ has to be understood whenever
$\alpha=0$ is chosen.) Let us provide a positive example of
concave-$\alpha$ functions in our context.

\begin{example}[Cumulative distribution functions]
Let $F \fc{\Re_+}{[0,1]}$ be the cumulative distribution function
of a $\chi$-random variable with $m$ degrees of freedom. Then for
any $\alpha$, $F$ is concave-$\alpha$ on the interval $I = (0, (m
- \alpha)^{\frac \alpha 2}]$ if $\alpha < 0$, $I =
[\ln(m)/2,\infty)$ if $\alpha = 0$ and $I = [(m - \alpha)^{\frac
\alpha 2}, \infty)$ if $\alpha \in (0,1]$.
This can be established by direct computation or as a result of
\cite[Lemma 3.1]{vanAckooij_Malick_2017}. Further positive
examples can be found in \cite[Table 1]{Henrion_Strugarek_2008}.\hfill\qedsymbol
\end{example}

As we will shortly see, $\invconcavityG$ and concavity-$\alpha$ of
cumulative distribution functions can be conveniently related to
specific properties of their density functions (provided they
exist). To this end, we introduce the following concept.

\begin{definition}[$G$-decreasing functions]
Let $G \fc{\Re_+}{\Re_+}$ be a strictly decreasing
(resp.\;increasing) continuously differentiable map with finitely
many critical points. A mapping $f$ is said to be $G$-decreasing
(resp.\;G-increasing) if there exists
$$t^\ast_{G} > [\max\col{t \; : G'(t)=0}]_+$$ such that the ratio
$r(t) := \frac{f(t)}{G'(t)}$ is strictly decreasing (increasing)
on the set $t \geq t^\ast_{G}$. Here $[t]_+ := \max\col{t,0}$ is
the positive part. \label{def:ratio}
\end{definition}

The instantiation of this definition, related to the family of
mappings $\col{\galpha}_\alpha$ was already introduced in
\cite{Henrion_Strugarek_2008} under the notion of
$\alpha$-decreasing functions (only considering the situation
$\alpha < 0$). We now provide a key result relating
$\invconcavityG$ of a given distribution function with its density
being $G$-decreasing.

\begin{proposition}[$\invconG$ cumulative distribution functions and $G$-decreasing densities]
Let $F \fc{\Re}{[0,1]}$ be the cumulative distribution function of
a random variable with associated (continuously differentiable)
density function\;$f$. Consider the statements:
\begin{enumerate}
\item the density $f$ is $G$-decreasing (see Definition
\ref{def:ratio}) with associated parameter $t^*$ ;
\item the mapping $F$ is $\invconG$ on the interval $I =
(0,G(t^\ast)]$ if $G$ is strictly decreasing, and on $I =
[G(t^\ast),\infty)$ if $G$ is strictly increasing, i.e., $z \maps
F(\reciprocal{G}(z))$ is concave on $I$ ;
%
\end{enumerate}
Then 1. implies 2. and if moreover $G$ is twice continuously
differentiable, then 2. also implies 1.. \label{prop:density}
\end{proposition}

\begin{proof}
\textbf{1. $\Rightarrow$ 2.}. We note that the proof of this
implication follows closely the proof of \cite[Lemma
3.1]{Henrion_Strugarek_2008} as well as \cite[Lemma
4]{Zadeh_Khorram_2012}. Let $G \fc{\Re_+}{\Re_+}$ be a strictly
increasing (decreasing) map and $t^\ast$ be such $t \maps
\frac{f(t)}{G'(t)}$ is strictly decreasing (increasing) on the set
$t \geq t^\ast$.

Let us begin by considering the situation wherein $G$ is strictly
decreasing. Then for $z \in (0, G(t^*))$, we have
$\reciprocal{G}(z) \geq t^\ast$. Now, the map $z \maps
F(\reciprocal{G}(z))$ that we will call $\chi$ can be written:
\begin{align*}
\chi(z) &= \int_{-\infty}^{\reciprocal{G}(z)} f(s) ds = F(t^\ast)
+
\int_{t^\ast}^{\reciprocal{G}(z)} f(s)ds \\
&= F(t^\ast) + \int_{G(t^\ast)}^{z}
\frac{f(\reciprocal{G}(u))}{G'(\reciprocal{G}(u))}du = F(t^\ast) -
\int_{z}^{G(t^\ast)}
\frac{f(\reciprocal{G}(u))}{G'(\reciprocal{G}(u))}du,
\end{align*}
where we have carried out the substitution $u = G(s)$. The ratio
appearing in the integral is a continuous map, making $\chi$
(continuously) differentiable. Moreover,
\begin{equation*}
\chi'(z) = \frac{f(\reciprocal{G}(z))}{G'(\reciprocal{G}(z))},
\end{equation*}
which together with $z \in (0,G(t^\ast))$ implies
$\reciprocal{G}(z) \geq t^\ast$ so that $\chi'$ is strictly
decreasing. As a consequence $\chi$ is indeed concave.

Let us now consider the case wherein $G$ is strictly increasing.
Then for $z \in (G(t^\ast),\infty)$ it also holds  that
$\reciprocal{G}(z) \geq t^\ast$. We can write $\chi$ as
\begin{equation*}
\chi(z) = \int_{-\infty}^{\reciprocal{G}(z)} f(s) ds = F(t^\ast) +
\int_{t^\ast}^{\reciprocal{G}(z)} f(s)ds = F(t^\ast) +
\int_{G(t^\ast)}^{z}
\frac{f(\reciprocal{G}(u))}{G'(\reciprocal{G}(u))}du
\end{equation*}
where we have carried out the substitution $u = G(s)$. Now
\begin{equation*} \chi'(z) =
\frac{f(\reciprocal{G}(z))}{G'(\reciprocal{G}(z))},
\end{equation*}
which together with $z \in (G(t^\ast),\infty)$ implies
$\reciprocal{G}(z) \geq t^\ast$ so that $\chi'$ is strictly
decreasing. As a result, $\chi$ is concave.

\noindent \textbf{2. $\Rightarrow$ 1.} Let us assume to begin with
that $F$ is $\invconG$, on the interval $[G(t^*), \infty)$ for a
strictly increasing map $G$ and define $\chi(z) = F\circ
\reciprocal{G}(z)$. We first note that $G$ is strictly increasing
and (continuously differentiable) and hence by the classic inverse
function Theorem (e.g., \cite[Theorem
1A.1]{Dontchev_Rockafellar_2014}, $\reciprocal{G}$ is also
continuously differentiable and the identity $(\reciprocal{G})'(x)
= \frac{1}{G'\circ \reciprocal{G}(x)}$ holds.

Now, by assumption, $\chi$ is concave on $[G(t^*), \infty)$, and
for any $x \in [G(t^*), \infty)$, we have
\begin{align*}
    \chi'(x) &= (\reciprocal{G})'(x).(f\circ G)(x) \\
    \chi''(x) &= (\reciprocal{G})''(x).(f\circ G)(x) + (\reciprocal{G})'(x)^2.(f'\circ
    G)(x).
\end{align*}
We can rewrite the second derivative as follows:
\begin{equation*}
\chi''(x) = \frac{1}{(G'\circ \reciprocal{G})(x)} \left ( -
\frac{G''\circ \reciprocal{G}(x)}{(G'\circ \reciprocal{G}(x))^2}
f\circ \reciprocal{G}(x) + \frac{f'\circ \reciprocal{G}}{G'\circ
\reciprocal{G}(x)} \right ),
\end{equation*}
where we have used the identity $(\reciprocal{G})''(x) = -
\frac{G''\circ \reciprocal{G}(x)}{(G'\circ \reciprocal{G}(x))^3}$
resulting from differentiating twice in the identity
$G(\reciprocal{G}(x)) = x$ holding locally at any $x
> 0$ and thus in particular at any $x \in [G(t^*), \infty)$ since
$G(t^*) > 0$. Since $G$ is strictly increasing, so is
$\reciprocal{G}$ and consequently $G'\circ \reciprocal{G}(x) > 0$.
Concavity of $\chi$ implies in turn that $\chi''(x) \leq 0$, which
can be equivalently stated as:
\begin{equation}
\frac{-G''(t)}{(G'(t))^2}f(t) + \frac{f'(t)}{G'(t)} \leq 0,
\label{eq:key}
\end{equation}
where $t = \reciprocal{G}(x)$ and $x \geq G(t^\ast)$ if and only
if $t \geq t^\ast$. By defining $\psi : t \mapsto
\frac{f(t)}{G'(t)}$ and differentiating once, we obtain:
\begin{equation*}
\psi'(t) = \frac{f'(t)}{G'(t)} + f(t) (-1) \frac{G''(t)}{G'^2(t)}
=  - \frac{G''(t)}{G'^2(t)} f(t) + \frac{f'(t)}{G'(t)}.
\end{equation*}
Now by \eqref{eq:key} it follows that $\psi'(t) \leq 0$ for all $t
\geq t^\ast$ and hence by Definition \ref{def:ratio}, $f$ is
$G$-decreasing.

This situation wherein $G$ is strictly decreasing follows upon
observing that $\reciprocal{G}$ is also strictly decreasing and
that consequently $G'\circ \reciprocal{G}(x) < 0$ holds.
Hence, concavity of $\chi$ on the set $(0, G(t^\ast)]$, implies
$\psi'(t) \geq 0$ as was to be shown.
\end{proof}

When applying the previous result with the family
$\col{\galpha}_{\alpha < 1}$ of \eqref{eq:galpha}, we obtain the
following corollary.

\begin{corollary}[Characterisation in the case of $\galpha$]\label{cor:hs}
Let $\alpha \in (-\infty,1)$ be given and $F: \mathbb{R}
\rightarrow [0,1]$ be the cumulative distribution function of a
random variable with continuously differentiable density $f$. Then
we have the following equivalence
\begin{enumerate}
\item $f$ is $\galpha$-decreasing (i.e., $f$ is
$(1-\alpha)$-decreasing in the sense of Definition 2.2 in
\cite{Henrion_Strugarek_2008})
\item $F$ is concave-$\galpha$ (i.e., $F$ is
$\alpha$-revealed-concave in the sense of  Definition 3.1 in
\cite{vanAckooij_Malick_2017}).
\end{enumerate}
\end{corollary}

The implication \textit{1.\;$\Rightarrow$ 2.} of Corollary
\ref{cor:hs} for $\alpha < 0$ was already known and corresponds to
Lemma 3.1 in \cite{Henrion_Strugarek_2008}. However both the
extension to $\alpha \in [0,1)$ and the reverse implication are
novel. Especially the latter shows that, in principle, there is no
loss of generality in studying the properties of the density
instead of the cumulative distribution function $F$.%

As already mentioned, \cite[Table\,1]{Henrion_Strugarek_2008}
contains a large choice of usual distribution functions (normal,
exponential, Weibull, gamma, chi, Maxwell, etc...) with a
$(1-\alpha)$-decreasing density function for all $\alpha < 0$ and
an analytic expression for the parameter $t^*$ indicated in
Definition \ref{def:ratio}. Although these results may give the
impression that all cumulative distribution functions are
$\invcon{\galpha}$, this is not true as the following example
shows.

\begin{example}[Not $\invcon{\galpha}$ distribution function]
Let $f \fc{\Re_+}{\Re}$ be defined as $f(t) =
(\frac{\sin^2(t)}{t^2})\frac{2}{\pi}$. Let us first verify that
$f$ does indeed integrate to $1$. This follows recalling the
identity $\sin^2(t) = (1-\cos(2t))/2$ and by using integration by
parts, as well as by recalling that $\lim_{t \rightarrow \infty}
\frac{\sin(s)}{s} ds= \frac{\pi}{2}$.

Now, should $f$ be $\galpha$-decreasing for some $\alpha < 0$,
then it must hold by \eqref{eq:key} that
\begin{equation*}
-\frac{\alpha(\alpha-1)}{\alpha^2}t^{-\alpha}f(t) +
\frac{1}{\alpha}t^{1-\alpha}f'(t) \leq 0,
\end{equation*}
for $t \geq t^\ast$ for some $t^*$. Yet the previous inequality is
equivalent with $(1-\alpha)f(t) + tf'(t) \geq 0$ for $t \geq
t^\ast$.

Let us verify that this can not hold. Note that
\begin{align*}
f'(t) & = \left(\frac{-2 \sin^2(t)}{t^{3}} +
\frac{2\sin(t)\cos(t)}{t^2}\right) \frac{2}{\pi}.
\end{align*}
We verify the negativeness of expression $g(t) := (1-\alpha)f(t) +
tf'(t)$, after algebraic manipulations:
\begin{align*}
g(t) = \frac{2}{\pi} \frac{\sin(t)}{t}
\left((-1-\alpha)\frac{\sin(t)}{t} + 2\cos(t) \right)
\end{align*}

Case $-1 < \alpha < 0$: Choosing the points $t_n = \frac{\pi}{2} +
2\pi n$, for integers $n>0$, we note that $\cos(t_n)=0$,
$\frac{\sin(t_n)}{t_n}>0$ and $(-1-\alpha)\frac{\sin(t_n)}{t_n} <
0 $. Hence, there is always $t_n$ such that $g(t_n) < 0$.

\noindent
Case $\alpha \leq -1$: Choosing the points $t_n = \frac{3\pi}{4} +
2\pi n$, for integers $n>0$ sufficient large such that $t_n >
-1-\alpha$. By noting that $\cos(t_n) = -0.71$ and $\sin(t_n) =
0.71$, it is easy to verify that $(-1-\alpha)\frac{\sin(t_n)}{t_n}
+  2\cos(t_n) < 0$. Again, consequently, there is always $t_n$
such that $g(t_n) < 0$.

Would the requested $t^\ast$ exist, we must have for some $n$
sufficiently large that $t_n > t^\ast$ and consequently condition
\eqref{eq:key} must hold in particular at $t=t_n$. Yet, we have
established that it can not. Hence the density function is not
$\galpha$-decreasing for any $\alpha < 0$.\hfill\qedsymbol
\end{example}

\begin{example}[$\invconcavityG$ with $G \neq \galpha$]
\label{ex:unusualG} Let us come back to Example\;\ref{ex:tbigalp}
and the map $G \fc{\Re_+}{\Re_+}$ defined as $G(x) =
\exp(-\ln(x)^{\frac 13})$. Then, $\Phi$ the cumulative
distribution function of a standard normal Gaussian random
variable is $\invconG$ on the set $(0,G(t^*)]$, with $t^\ast =
1.86$.

This follows from Proposition \ref{prop:density} as soon as the
ratio $r(x) = \frac{1}{\sqrt{2\pi}G(x)}\exp(-\frac 12
x^2)(-3x\ln(x)^{\frac 32})$ is strictly increasing. This, in turn,
can be asserted if the function $f(x) := \exp(-\frac 12 x^2)(-3x\ln(x)^{\frac
32})$ is strictly increasing on the set $x \geq t^*$. In order to show this, compute the derivative: $f'(x) = 3
\exp(-\frac 12 x^2)\ln(x)^{\frac 12}((x^2-1)\ln(x)-\frac 32)$.
Observing that for $x \in (1,\infty)$, $(3 \exp(-\frac 12
x^2)\ln(x)^{\frac 12}) > 0$, the sign of  $f'(x)$ depends on the
term: $f_2(x) = (x^2-1)\ln(x)-\frac 32$. The latter has
derivative: $f'_2(x) = 2x\ln(x) + \frac{(x^2 -1)}{x}$.
For $x > 1$, $f'_2(x) > 0$ ($f_2(x)$ is strictly increasing for
$x>1$). We find that $f_2(1.8) = -0.1834$ and $f_2(1.9) = 0.1752$,
so there is a root of $f_2$ in the interval $[1.8, 1.9]$. So, for
$x \in [1.9,\infty)$, $f_2(x) > 0$ and hence $f'(x)
> 0$, which in turn implies that $f(x)$ is strictly increasing. Numerically solving
$f_2(x) =0$ in\;$x$, we find the solution $x^*=1.8528$.\hfill\qedsymbol
\end{example}

%
%
%
%

\section{Interplay of generalized concavity and ``convexity" of chance constraints}
\label{sec:interplay}

In this section, we establish convexity results for feasible sets
of probabilistic constraints by employing the
set of tools of generalized concavity. Our analysis considers two
special structures for \eqref{eq:prb}. The first situation,
analyzed in Section \ref{subsec:nonlin}, refers to the general
case wherein $g$ is non-linear and relatively arbitrary, but the
random vector $\xi$ is assumed to follow a multi-variate
elliptically symmetric distribution. The second situation,
analyzed in Section \ref{subsec:copulae}, is when $g$ is
separable, which boils\footnote{A general separable $g$ would be
of the form $g(x,z) = \psi(z) - h(x)$. In our case, recalling that
$z$ will be substituted out for the random vector $\xi$, there is
no loss of generality in assuming $\psi(z)=z$ since the general
case reduces to it when taking $\tilde \xi = \psi(\xi)$ to be the
underlying random vector that we study.} down to considering
$g(x,z) = z - h(x)$. The random vector $\xi$ can be relatively
arbitrary in so much that it can have nearly arbitrary marginal
distributions and the (joint) dependency structure is pinned down
by the choice of a copula.

\subsection{Non-linear couplings of decisions vectors and elliptically distributed random vectors}
\label{subsec:nonlin}

In this section we consider the situation of \eqref{eq:prb}
wherein the map $g \fc{X \times \Re^m}{\Re^k}$ is convex in the
first argument and continuous as a function of both arguments. We
also assume that the random vector $\xi$ taking values in
$\mathbb{R}^m$ is elliptically symmetrically distributed with mean
$\mu$, covariance-like matrix $\Sigma$ and generator $\theta :
\mathbb{R}_+ \to \mathbb{R}_+$, which is denoted by  $\xi \sim
\mathcal{E}(\mu,\Sigma,\theta)$ if and only if its density $f_\xi
: \mathbb{R}^m \to \mathbb{R}_+$ is given by
    \begin{align}
    f_\xi (z) = \big( \det \Sigma \big)^{-1/2} \theta \bigg(  (z- \mu)^{\top} \Sigma^{-1} (z- \mu)
    \bigg),
    \end{align}
where the generator function $\theta \fc{\Re_+}{\Re_+}$ must
satisfy
\begin{equation*}
\int_0^\infty t^{\frac m2}\theta(t) dt < \infty.
\end{equation*}
We consider $L$ as the matrix arising from the Choleski
decomposition of $\Sigma$, i.e., $\Sigma = LL^\mtr$, it can be
shown that $\xi$ admits a representation as
\begin{equation}\label{eq:radial}
\xi = \mu + \mathcal{R} L \zeta
\end{equation}
where $\zeta$ has a uniform distribution over the Euclidean
$m$-dimensional unit sphere $S^{m-1}:=\{ z\in \mathbb{R}^m :
\sum_{i=1}^m z_i^2=1 \}$ and $\mathcal{R}$ possesses a density,
which is given by
\begin{equation}\label{densityR}
f_{\mathcal{R}}(r):=\frac{2 \pi^\frac{m}{2} }{\Gamma(\frac{m}{2})}
r^{m-1} \theta(r^2),
\end{equation}
with $\Gamma$ is the usual gamma-function.

The family of elliptically distributed random vectors includes many
classical families (see e.g.\,\cite{Fang_Kotz_Ng_1990} and
\cite{Landsman_Valdez_2013}): for instance, Gaussian random
vectors and Student random vectors (with $\nu$ degrees of freedom)
are elliptical with the respective generators
\begin{equation*}
\theta^{\rm Gauss} (t)= \exp(-t/2)/(2\pi)^{m/2}
\qquad\text{and}\qquad \theta^{\rm Student} (t)=
\frac{\Gamma\big(\frac{m+\nu}{2}\big)}{\Gamma\big(\frac{\nu}{2}\big)}(\pi\nu)^{-m/2}\big(1+\frac{t}{\nu}\big)^{-\frac{m+\nu}{2}}.
\end{equation*}

The advantage of the spherical radial decomposition
\eqref{eq:radial} is that it allows one to derive
the following attractive form for $\phi$ (see e.g.\;Theorem\,2.1 of\;\cite{vanAckooij_Malick_2017}): if $x \in X$ is such that
\begin{enumerate}
\item $g(x,\mu) \leq 0$ (recall that $\mu = \esp{\xi}$)
\item for any $z \in \Re^m$ such that $g(x,z) \leq 0$, we have
\begin{equation*}
g(x, \lambda \mu + (1-\lambda)z) \leq 0 \foral \lambda \in [0,1],
\end{equation*}
\end{enumerate}
then $\phi$ defined in \eqref{eq:prb} can be written as
\begin{equation}
\phi(x) = \int_{v\in \mathbb{S}^{m-1}} F_{\scr{R}}(\rho(x,v))
d\mu_{\zeta}(v) \label{eq:phirep}
\end{equation}
where $F_{\scr{R}}$ is the cumulative distribution function of
$\scr{R}$, $\mu_\zeta$ is the law of uniform distribution on the $m$-dimensional euclidian sphere 
$\mathbb{S}^{m-1}$, and $\rho \fc{X \times \mathbb{S}^{m-1}}{\mathbb{R}_+
\cup \col{\infty}}$ is the continuous\footnote{The continuity of $\rho$ is well-known (see e.g.\;Lemma 3.4 in \cite{van2020generalized}); it implies in particular that $\rho$ is mesurable. Note also that $\rho$ is quasi-concave (see Lemma 3.2 in\;\cite{vanAckooij_Malick_2017}). We will need here stronger notions of $G$-concavity to establish our results.}  mapping defined by
\begin{equation}
\rho(x,v) = \left\{
\begin{aligned}
&\sup_{t\geq 0}& &t\\
&s.t. &g(x,\mu + &tLv) \leq 0. \\
\end{aligned}
\right. \label{eq:rhodef}
\end{equation}

Note that if
for each $v$ the map $\rho(\cdot,v)$ is $F_{\scr{R}}$-concave,
then due to the linearity of the integral such a property would
carry over immediately to $\phi$. It is clear that such a request
could not hold without restrictions since generally a probability
function can not be ``concave''. Indeed, it is a bounded function
by $0$ and $1$ and usually increasing along a certain ``path''.
For this reason the analysis is non-trivial. The second difficulty
is in conveying these desired properties from $p$ alone. Let us
provide a precise statement.

\begin{theorem}[Convexity of probability functions]
Let $C$ be a convex subset of $X$. Assume that, for any $v \in
\sph$, there exists a continuous function $G_v : \rho(C \times
\col{v}) \rightarrow \mathbb{R_+}$ such that
\begin{itemize}
\item $G_v$ is strictly monotonic on $\rho(C \times \{ v \})$,
\item $x \mapsto \rho(x,v)$ is $G_v$-concave on $C$ and
continuous,
\item $F_{\scr{R}}$ is $\invcon{G_v}$ on $(G_v \circ \rho)(C
\times \{ v \})$,
\end{itemize}
where $\rho$ is defined as in \eqref{eq:rhodef}. Then $\phi
\fc{X}{[0,1]}$ defined in \eqref{eq:prb} is concave on $C$.
\label{thm:phiccv}
\end{theorem}

\begin{proof}
Let us first establish that for any fixed $v \in \sph$ that $x
\mapsto F_{\scr{R}}\circ \rho(x,v)$ is concave on $C$. To this end
pick $x_1, x_2 \in C$ and $\lambda \in [0,1]$ arbitrarily and
consider $x^\lambda = \lambda x_1 + (1-\lambda)x_2$. Then by
$G_{v}$-concavity of $\rho(\cdot,v)$ it follows:
\begin{equation*}
\rho(x^\lambda,v) \geq \reciprocal{G_{v}}(\lambda G_v(\rho(x_1,v))
+ (1-\lambda) G_v(\rho(x_2,v)) ).
\end{equation*}
Now since $F_{\scr{R}}$ is increasing as a distribution function
it follows too that
\begin{equation*}
F_{\scr{R}}(\rho(x^\lambda,v)) \geq F_{\scr{R}}(
\reciprocal{G_{v}}(\lambda G_v(\rho(x_1,v)) + (1-\lambda)
G_v(\rho(x_2,v)) ) ).
\end{equation*}
Now consider the set $I_v \subseteq \Re$ defined as $I_v =
G_v(\rho(C,v))$.

Our claim is that $I_v$ is an interval, i.e., is convex. To this
end, we recall that the continuous image of a connected set is
connected and hence $\rho(C,v)$ is a connected set (recall that
$C$ is a convex set). By applying the argument a second time,
since $G_v$ is continuous, it follows that $I_v$ is connected. Now
since $I_v$ is a subset of $\dl{R}$, it is connected if and only
if it is convex if and only if it is an interval.

Now since $I_v$ is an interval and hence convex, $\lambda
G_v(\rho(x_1,v)) + (1-\lambda) G_v(\rho(x_2,v)) ) \in I_v$ and we
may apply $\invconcavity{G_v}$ of $F_{\scr{R}}$ to pursue our
development as follows:
\begin{equation*}
F_{\scr{R}}(\rho(x^\lambda,v)) \geq (F_{\scr{R}} \circ
\reciprocal{G_{v}})(\lambda G_v(\rho(x_1,v)) + (1-\lambda)
G_v(\rho(x_2,v)) ) ) \geq \lambda F_{\scr{R}}(\rho(x_1,v)) +
(1-\lambda) F_{\scr{R}}(\rho(x_2,v)),
\end{equation*}
which is what was to be shown. Now by linearity of integrals, we
get
\begin{align*}
\phi(x^\lambda) &= \int_{v\in \mathbb{S}^{m-1}}
F_{\scr{R}}(\rho(x^\lambda,v)) d\mu_{\zeta}(v) \geq  \lambda
\int_{v\in \mathbb{S}^{m-1}} F_{\scr{R}}(\rho(x_1,v))
d\mu_{\zeta}(v) + (1-\lambda)\int_{v\in \mathbb{S}^{m-1}}
F_{\scr{R}}(\rho(x_2,v)) d\mu_{\zeta}(v) \\
&= \lambda\phi(x_1) + (1-\lambda)\phi(x_2),
\end{align*}
thus concluding the proof.
\end{proof}

Although Theorem \ref{thm:phiccv} allows us to establish concavity
of $\phi$ on a certain given convex set $C$, an important
additional difficulty is how to entail that $M(p) \subseteq C$
holds for $p$ large enough. Then one can immediately deduce the
convexity of $M(p)$ from concavity of $\phi$. A convenient
situation is one when, for all $x$, the set $\mathfrak{M}(x) :=
\{z \in \mathbb{R}^{m} : g(x,z) \leq 0\}$ is convex in $\Re^m$.
For this it would be sufficient to request that $g$ is convex
respectively in $x$ and in $z$ (but not necessarily jointly). We
can however generalize to the situation wherein $\mathfrak{M}(x)$
is star-shaped with respect to $\mu$ if the convex hull of the
latter sets does ``not distort" measurement of length of lines
segments $\col{r \geq 0 \; : \mu + rLv \in \frk{M}(x)}$ moving
through it. In order to make a precise statement, we introduce the
following map: $\rho^{\tt co} \fc{X \times \sph}{\Re_+ \cup
\col{-\infty, \infty}}$:
\begin{equation}
\rho^{\tt co}(x,v) = \left\{
      \begin{aligned}
        &\sup_{t\geq 0}& &t\\
        &s.t. &\mu + &tLv \in \hull(\mathfrak{M}(x))\\
      \end{aligned}
    \right. \label{eq:rhocodef}
\end{equation}
where $\hull(\mathfrak{M}(x))$ denotes the convex hull of
$\mathfrak{M}(x)$.

\begin{theorem}[Eventual convexity of elliptical chance constraints]
\label{thm:evencvx} Let $C$ be a given convex subset of $X$. In
addition to the framework of this section, assume that
\begin{enumerate}
\item there exists a $t^* > 0$, such that $\{x \in X \; :\rho(x,v)
\geq t^*, \forall v \in \mathbb{S}^{m-1} \} \subseteq C$\;;
\item For any $v \in \sph$, there exists a continuous function
$G_v : \mathbb{R} \rightarrow \mathbb{R}$ as in Theorem\;\ref{thm:phiccv};
%
\item There exists $p_0 \in$ $[\frac{1}{2},1]$ and $\delta^{\tt
nd}
> 0$ such that
\begin{equation*}
\delta^{\tt nd}\rho(x,v) \geq \rho^{\tt co}(x,v),
\end{equation*}
for all $x \in M(p_0)$ and all $v \in \sph$, where $\rho$ and $\rho^{\tt co}$ are defined respectively in
\eqref{eq:rhodef} and \eqref{eq:rhocodef}.
\end{enumerate}

Then for any $q \in (0,\frac{1}{2})$ and any $p \geq \max(p_0,
p(t^\ast,q))$ with
\begin{equation}
p(t^\ast,q) = \left(\frac{1}{2} - q\right)
F_{\scr{R}}\left(\frac{\delta^{\tt nd}t^*}{\delta(q)}\right) +
\frac{1}{2} + q,
\end{equation}
the set $M(p)$ defined in \eqref{eq:Mp} is convex. Here
$\delta(q)$ is the unique solution (in $\delta$) to the equation
\begin{equation*}
\frk{B}_i\left(\frac{m-1}{2}, \frac 12,
\sin^2(\arccos(\delta))\right) =
(1-2q)\,\frk{B}_c\left(\frac{m-1}{2}, \frac 12\right),
\end{equation*}
where $\frk{B}_i$ (resp.\;$\frk{B}_c$) refers to the incomplete
(resp.\;complete) Beta function.
\end{theorem}

\begin{proof}
We follow here closely the demonstration of the Theorem 4.1 from
\cite{vanAckooij_Malick_2017}. Let $p \in [p_0,1]$ be given and
take any $x \in M(p)$. We have, for such $x$ that $\frac{1}{2} < p
\leq \prb[g(x,\xi) \leq 0] \leq \prb[\xi \in
Co(\mathfrak{M}(x))]$. Then corollary 2.1 from
\cite{vanAckooij_Malick_2017} gives that $\mu \in
int(Co(\mathfrak{M}(x)))$.

Let us now pick an arbitrary but fixed $v \in \dom{\rho(x, .)}$.
Note that, by assumption, there is $\delta^{\tt nd} > 0$ such that
$\delta^{\tt nd}\rho(x,v) \geq \rho^{\tt co}(x,v)$, i.e., $v \in
\dom{\rho^{\tt co}(x, .)}$ as well. Hence, $\mu + \rho^{\tt
co}(x,v)Lv$ belongs to the boundary of $\hull(\mathfrak{M}(x))$.

Therefore, we can separate $\mu + \rho^{\tt co}(x,v)Lv$ from the
convex set $\hull(\mathfrak{M}(x))$, so that there exists a
non-zero $s \in \Re^m$ such that for all $z \in
\hull(\mathfrak{M}(x))$,
\begin{equation*}
s^\mtr z \leq s^\mtr(\mu + \rho^{\tt co}(x,v) Lv) \leq s^\mtr(\mu
+ \delta^{\tt nd}\rho(x,v) Lv).
\end{equation*}
Now define $c \in \mathbb{R}^m$ and $\gamma > 0$ as follows:
\begin{equation*}
c := \frac{s}{\norm{L^\mtr s}} \; \mbox{and}\; \gamma = c^\mtr(\mu
+ \delta^{\tt nd} \rho(x,v) Lv),
\end{equation*}
where we recall the $\norm{L^\mtr s} > 0$ since $L$ is regular and
$s \neq 0$. It now follows by construction that,
\begin{equation*}
\mathfrak{M}(x) \subset \hull(\mathfrak{M}(x)) \subset \{z \in
\mathbb{R}^m \; : c^\mtr z \leq \gamma \}.
\end{equation*}

In particular this entails $\prb[g(x, \xi) \leq 0] \leq
\prb[c^\mtr \xi \leq \gamma]$. We can employ Theorem 2.2 of
\cite{vanAckooij_Malick_2017} to get the estimate
\begin{align*}
p &~\leq~ \prb[g(x, \xi) \leq 0] ~\leq~ \prb[c^\mtr \xi \leq \gamma] ~\leq~ \left(\frac{1}{2} - q\right) F_{\scr{R}}\left(\frac{\delta^{\tt nd} \rho(x,v)\frac{s^\mtr Lv}{\norm{L^\mtr s}}}{\delta(q)}\right) + q + \frac{1}{2} \\
  &\leq \left(\frac{1}{2} - q\right) F_{\scr{R}}\left(\frac{\delta^{\tt nd} \rho(x,v)}{\delta(q)}\right) + q +
  \frac{1}{2},
\end{align*}
for any $q \in (0, \frac 12)$ and associated $\delta(q) > 0$,
where we have used the Cauchy-Schwartz inequality and the
monotonicity of $F_{\scr{R}}$. Since $F_{\scr{R}}$ is increasing
and $(\frac{1}{2}-q) > 0$ as well as $p \geq \max(p_0, p(t^\ast,
q))$, we derive $\rho(x,v) \geq t^*$. Moreover obviously,
$\rho(x,v) \geq t^\ast$ for $v \notin \dom{\rho(x,.)}$. Hence,
$M(p) \subseteq \{x \in X \; : \rho(x,v) \geq t^*, \forall v \in
\mathbb{S}^{m-1} \} \subseteq C$. Now, we can apply Theorem
\ref{thm:phiccv} to establish that $\phi$ is concave on $C$ and
therefore $M(p)$ must be convex.
\end{proof}

\begin{remark}[Abstract theorem at work]
Though looking abstract, the conditions of the theorem are often
present in practice. They can for example be ensured whenever the
following conditions hold. We provide examples in
Section\;\ref{sec:examples}.
\begin{itemize}
\item There exist some $\alpha \in \Re$ such that for each $v \in
\sph$, the map $x \maps \rho(x,v)$ is $\alpha$-concave (see
\cite[Proposition 5.1]{vanAckooij_Malick_2017} for an exemple) and
the radial distribution function $F_{\scr{R}}$ is
concave-$\alpha$. Here prominent examples are the chi
distribution, the Fisher-Snedecor distribution, etc.
\item To get a suitable $t^\ast$ associated with $F_{\scr{R}}$,
one can use the concavity-$\alpha$ of $F_{\scr{R}}$. For the
example of the chi-distribution, we get $t^\ast = \sqrt{m -
\alpha}$ and $C = \col{x \in X \; : \rho(x,v) \geq \sqrt{m -
\alpha} \ \text{for all } v \in \sph}$.
\item The request of item 3 holds whenever the set $\frk{M}(x)$ is
convex for all $x \in C$ and in that case $\delta^{\tt nd}=1$.
Such convexity can be ensured whenever the map $z \maps g(x,z)$ is
quasi-convex for each $x \in X$, and in that case $p_0 = \frac 12$
can be taken.

\item Note finally that continuity of the map $\rho$ can be
ensured under fairly general conditions, e.g., convexity of $g$ in
the second argument and $g(x,\mu) < 0$ together with continuous
differentiability of $g$ immediately entail continuity of $x \maps
\rho(x,v)$ (even continuity in both arguments). See for instance
\cite{Hantoute_Henrion_Perez-Aros_2017}. Continuity can also be
ensured under less restrictive conditions. For instance whenever,
the map $g$ is continuous at any $\bar x, \bar v, \bar r$ such
that $g(\bar x, \bar r L \bar v) = 0$, neighbourhoods $U, V, W$ of
$\bar x, \bar v, \bar r$ respectively can be identified such that
for all $(x, v) \in U \times V$, the map $r \maps g(x,rLv)$ is
monotonic. This in turn is related to uniqueness of solutions of
perturbed systems $g(x,rLv) = t$ for $t$ sufficiently small and
not very restrictive. To ensure continuity of $\rho$ from this
assumption, we can just use a general version of the implicit
function theorem; see \cite[Theorem
1H.3]{Dontchev_Rockafellar_2014} as well as
\cite{Jittorntrum_1978,Kumagai_1980} which are older statements of
such a result.\end{itemize}
\end{remark}

The computed threshold $p^*$ is thus valid for a large class of
nonlinear functions $g$. Should it be conservative,  refinements
might be obtained when studying specific structures; see e.g.,
\cite{Minoux_Zorgati_2016,vanAckooij_2016b} and forthcoming
Remark\;\ref{rem:improve}. Beyond being a general guarantee, the
threshold $p^*$ is thus an indication that lower thresholds could
be revealed for specific functions.

\subsection{Separable copul\ae\/-structured probabilistic constraints}
\label{subsec:copulae}

In this section we will consider the following form of
\eqref{eq:prb}:
\begin{equation*}
\prb[ \xi \leq h(x)] \geq p,
\end{equation*}
where $\xi$ taking values in $\Re^m$ is a random vector and $h : X
\rightarrow \mathbb{R}^m$ a given map. By employing Sklar's
Theorem \cite{Sklar_1959}, we may write \eqref{eq:prb} with the
special structure as above in the following form:
\begin{equation}
\mathbb{P} [ \xi \leq h(x)] = \C(F_1(h_1(x)),...,F_m(h_m(x))).
\label{eq:copCC}
\end{equation}
Here $F_1,...,F_m$ are the marginal distribution functions of the
random vector $\xi$ and $h_i \fc{X}{\Re}$, $i=1,...,m$ refers to
the components of the mapping $h$. Moreover $\C
\fc{[0,1]^m}{[0,1]}$ is a copula, i.e., a multi-variate
distribution function with uniform marginal distributions (e.g.,
\cite{Nelsen_2006} for further details). In what follows, we
require specific properties of $\C$ and $F_1,...,F_m$ and through these choices pin down the
multivariate distribution of\;$\xi$.

We introduce the analogue of $\invconG$ functions in this context
of copul\ae\/.
\begin{definition}[$\invconG$ copul\ae\/]\label{def:copula}
Let $\C \fc{[0,1]^m}{[0,1]}$ be a copula and $G
\fc{[0,1]^m}{\Re^m}$ a map such that the $i$th component $G_i$ is
continuous and strictly monotonic. The copula $\C$ is said to be
$\invconG$ on the product of intervals $I = \prod_{i=1}^m I_i$,
with $I_i \subseteq \Re$ if the map $$I \ni z \maps
\C(\reciprocal{G_1}(z_1),...,\reciprocal{G_m}(z_m))$$ is
quasi-concave.
\end{definition}

\begin{example}[Relation to other concave copul\ae\/]\label{ex:relation}
Taking $G_i=G_{\gamma}$ in the previous definition, a $\invconG$
copula corresponds to a $(-\infty)-\gamma$-concave copula in the
terminology of \cite{vanAckooij_2013}. For example, all
Archimedian copul\ae\/ are $\invcon{G_1}$ (see \cite[Theorem
3.3]{vanAckooij_Oliveira_2016}. Example \ref{ex:GasCopConc} below
provides an example of $\invconG$ copula which is necessarily
$(-\infty)-\gamma$-concave. Note also that, in the special case
$\gamma=0$, this notion is weaker than the notion of
logexp-concave
of \cite{Henrion_Strugarek_2011}). Prominent examples of such
copul\ae\/ are for instance the independent, maximum or Gumbel
copula. The Clayton copula is an example of $\invcon{G_{0}}$
copula which is not logexp-concave (see Lemma\;5.5
of\;\cite{vanAckooij_2013}).\hfill\qedsymbol
%
%
\end{example}

\begin{example}[Gaussian case]\label{ex:GasCopConc}
Let $R$ be a positive definite $m \times m$ correlation matrix and
$\Phi$ denote the standard normal distribution function. Then the
Gaussian copula $\C^R \fc{[0,1]^m}{[0,1]}$ is defined as
\begin{equation*}
\C^R(u) =
\Phi^R(\reciprocal{\Phi}(u_1),...,\reciprocal{\Phi}(u_m)),
\end{equation*}
where $\Phi^R$ is the multivariate Gaussian distribution function
related to correlation matrix $R$. Now observe that
$\reciprocal{\Phi} \fc{[0,1]}{\Re}$ is strictly increasing and
continuous. Now for any $z \in I = \Re^m$,
$\C^R(\Phi(z_1),...,\Phi(z_m)) =
\Phi^R(\reciprocal{\Phi}(\Phi(z_1)),...,\reciprocal{\Phi}(\Phi(z_m)))
=\Phi^R(z_1,...,z_m)$. Recalling that multivariate Gaussian
distribution functions are $0$-concave (see, e.g.,
\cite{Prekopa_1995}), it follows from Corollary \ref{cor:alpconc}
that $\C^R$ is $\invcon{\reciprocal{\Phi}}$. It was however not
clear whether or not $\C^R$ is $(-\infty)-\gamma$-concave; see the
extensive discussion in \cite[section 6]{vanAckooij_2013}).
Similar analysis can be carried out with Clayton Copula which can
be shown to be $\invconG$ with $G_i(z) = \ln{(z)^2}$.\hfill\qedsymbol
\end{example}

We can now provide the announced eventual convexity result.
\begin{theorem}[Eventual convexity of separable copulae-structured probabilistic constraints]\label{thm:convexity}
Let $h_i \fc{X}{\dl{R}}$ be continuous mappings and consider the
following identity:
\begin{equation}
\prb[\xi \leq h(x)] = \C(F_1(h_1(x)),...,F_m(h_m(x))),
\end{equation}
where $\C$ is a suitable copula and $F_i$ are the marginal
distribution functions of component $i$ of $\xi$, $i=1,...,m$.

Assume that we can find strictly monotonous mappings $G_i
\fc{\Re}{\Re}$, such that the functions $h_i$ are $G_i$-concave on
a given convex level set $C = \col{x \in X \; : h_i(x) \geq b_i,
i=1,...,m}$ of $X$ for appropriate parameters $b_i \in \Re$,
$i=1,...,m$.

Assume moreover that for $i=1,...,m$ continuous strictly monotonic
mappings $\hat G_i \fc{[0,1]}{\Re}$ can be identified such that
$\C$ is $\invcon{\hat G}$ (see Definition \ref{def:copula}) on the
set $I = \prod_{i=1}^m I_i$, where
\begin{itemize}
\item the interval $I_i = [\hat G_i(b_i),\infty)$ whenever $\hat
G_i$ is strictly increasing
\item the interval $I_i = (-\infty, \hat G_i(b_i)]$ whenever $\hat
G_i$ is strictly decreasing.
\end{itemize}
Finally, assume that the marginal distribution functions $F_i
\fc{\Re}{\Re}$ are $\hat G_i$-concave-$G_i^{-1}$ on
\begin{itemize}
\item the interval $[G_i(b_i),\infty)$ whenever $G_i$ is strictly
increasing
\item the interval $(0, G_i(b_i)]$ whenever $G_i$ is strictly
decreasing.
\end{itemize}

Then the set $M(p) := \col{x \in X \; : \prb[\xi \leq h(x)] \geq
p}$ is convex for all $p > p^*:= \max_{i=1,...,m}{F_i(b_i)}$.
Convexity can also be asserted for $p = p^*$ if each individual
distribution function $F_i$, $i=1,...,m$ is strictly increasing.
\end{theorem}

\begin{proof}
Pick any $p > p^*$, $x,y\in M(p)$, $\lambda \in [0,1]$ and $i \in
\col{1,...,m}$ arbitrarily. Define $x^{\lambda} := \lambda x +
(1-\lambda)y$. Since all copul\ae\/ are dominated by the
maximum-copula, we get:
\begin{equation}
F_i(h_i(x)) \geq \min_{j=1,...,m} F_j(h_j(x)) \geq
\C(F_1(h_1(x)),...,F_m(h_m(x))) \geq p > p^* \geq F_i(b_i).
\label{eq:pstar}
\end{equation}
Now the latter entails
\begin{equation}
h_i(x) \geq b_i. \label{eq:pstar2}
\end{equation}
and this in turn means $M(p) \subseteq C$. Estimate
\eqref{eq:pstar2} also holds whenever $p\geq p^*$ and $F_i$ is
strictly increasing for each $i=1,...,m$. A similar estimate is
obtained for $y$ clearly. Let us first remark that due to
\eqref{eq:pstar2}, we have $G_i(h_i(x)) \geq G_i(b_i)$, whenever
$G_i$ is strictly increasing and $G_i(h_i(x)) \leq G_i(b_i)$
otherwise.

Consequently, $\lambda G_i(h_i(x)) + (1-\lambda)G_i(h_i(y)) \leq
G_i(b_i)$ whenever $G_i$ is strictly decreasing, and $\lambda
G_i(h_i(x)) + (1-\lambda)G_i(h_i(y)) \geq G_i(b_i)$ otherwise. We
note that in either situation $\lambda G_i(h_i(x)) +
(1-\lambda)G_i(h_i(y))$ belongs to the interval associated with
$\hat G_i$-concavity-$G_i^{-1}$ of $F_i$. Now, since $F_i$ is
increasing we establish first that
\begin{equation}
F_i(h_i(x^{\lambda}))\geq F_i( \reciprocal{G_i}(\lambda
G_i(h_i(x)) + (1-\lambda)G_i(h_i(y))) ), \label{eq:aa}
\end{equation}
by $G_i$-concavity of $h_i$ on the set $C$. As argued above,
related to $\lambda G_i(h_i(x)) + (1-\lambda)G_i(h_i(y))$
belonging to an appropriate domain, we may pursue \eqref{eq:aa} by
invoking $\hat G_i$-concavity-$G_i^{-1}$ of $F_i$, and obtain
\begin{align}
F_i(h_i(x^{\lambda}))& \geq F_i( \reciprocal{G_i}(\lambda G_i(h_i(x)) + (1-\lambda)G_i(h_i(y))) ) \notag \\
&\geq \reciprocal{\hat G_i}( \lambda \hat G_i( F_i
(\reciprocal{G_i}(G_i(h_i(x))))) + (1-\lambda) \hat G_i( F_i
(\reciprocal{G_i}(G_i(h_i(x))))) ) \notag \\
&\geq \reciprocal{\hat G_i}( \lambda \hat G_i( F_i(h_i(x)) ) +
(1-\lambda) \hat G_i( F_i(h_i(y)) ) ) \label{eq:ab}
\end{align}
Since $i$ was fixed but arbitrary, the above equation
\eqref{eq:ab} holds for all $i=1,...,m$. A Copula is increasing in
its arguments, so we get in turn:
\begin{align*}
& \C(F_1(h_1(x^\lambda)),...,F_m(h_m(x^\lambda))) \geq \nonumber
\\
& \quad \C( \reciprocal{\hat G_1}( \lambda \hat G_1( F_1(h_1(x)) )
+ (1-\lambda) \hat G_1( F_1(h_1(y)) ) ), ..., \reciprocal{\hat
G_m}( \lambda \hat G_m( F_m(h_m(x)) ) + (1-\lambda) \hat G_m(
F_m(h_m(y)) ) ) ).
\end{align*}
Moreover arguing as before, for any $i=1,...,m$ it holds that
$\lambda \hat G_i( F_i(h_i(x)) ) + (1-\lambda) \hat G_i(
F_i(h_i(y)) ) \in I_i$. We can now employ $\invconcavity{\hat G}$
of the copula $\C$, to establish:
\begin{align*}
& \C(F_1(h_1(x^\lambda)),...,F_m(h_m(x^\lambda))) \geq  \\
& \quad
m_{-\infty}(\C(F_1(h_1(x)),...,F_m(h_m(x))),\C(F_1(h_1(y)),...,F_m(h_m(y))),\lambda)
\geq p,
\end{align*}
which is equivalent with
$\C(F_1(h_1(x^\lambda)),...,F_m(h_m(x^\lambda))) \geq p$, i.e.,
$x^{\lambda} \in M(p)$ as was to be shown.
\end{proof}

The conditions given in the above Theorem \ref{thm:convexity}
allow for quite some flexibility. We refer to \cite[section
5]{vanAckooij_2013} for examples. In section \ref{sec:excop} below
we provide other examples, not covered by prior results.


\section{Selected examples}
\label{sec:examples}

In this section, we provide several new examples for which
eventual convexity can be asserted with the extended framework
built upon $G$-concavity. Section \ref{subsec:quad} is concerned
with a situation wherein $g$ is non-linear and Section
\ref{sec:excop} wherein $g$ separable.

\subsection{Example with a quadratic probabilistic constraint}
\label{subsec:quad}

Consider the map $g \fc{\Re^n \times \Re^m}{\Re}$ defined as
\begin{equation}
g(x,z) = z^\mtr W(x) z + 2\sum_{i=1}^n a_i w_i^\mtr z + b
\label{eq:gquad}
\end{equation}
where $a_1,...,a_n$ are fixed constants, $w_1,...,w_n$ are fixed
vectors of $\mathbb{R}^m$, and $W\colon\mathbb{R}^n \rightarrow
S_n^+ $ is a ``convex" function, in the sense that for any $(x,y)
\in \mathbb{R}^n, \lambda \in$ [0,1], $\lambda W(x) + (1-\lambda)
W(y) - W(\lambda x + (1-\lambda) y)$ is positive (semi-)definite.
Note that, for all $z$, $g(\cdot, z)$ is then convex with respect
to its first variable. We assume furthermore that $b < 0$, which
ensures $g(x,0) < 0$. A concrete example of a mapping $W$
satisfying the above request is for instance $W(x) = \sum_{i=1}^n
x_i W_i$, with $W_i$ positive semi-definite and $x \geq 0$ as an
ambiant further restriction.
%

Let $\xi$ be an elliptically symmetrically distributed random
vector with given spherical-radial decomposition $\xi = \scr{R}L
\zeta$. Then, by fixing $(x,v) \in \mathbb{R}^n\times
\mathbb{S}^{m-1}$ and defining
\begin{align*}
h(x,v) & := (Lv)^\mtr W(x)(Lv)\quad \text{and}\quad \beta(v) :=
\sum_{i=1}^n a_i w_i^\mtr (Lv),
\end{align*}
we can explicitly identify the map $\rho$ of \eqref{eq:phirep},
being the unique solution to the equation $g(x,rLv)=0$. This
solution map is given by the expression:
\begin{equation}
    \rho(x,v) = \frac{-\beta(v) + \sqrt{\beta(v)^2 - b
    h(x,v)}}{h(x,v)}.
\end{equation}

Next, we identify a function $g_v : \mathbb{R}^+ \rightarrow
\mathbb{R}$ such that $\rho(.,v)$ is $g_v$-concave for any $v \in
\mathbb{S}^{m-1}$.

\begin{lemma}
For any fixed $v \in \mathbb{S}^{m-1}$ define the map $g_v
\fc{(0,\infty)}{\Re_-}$ as $g_v(t) = -(\frac{b}{t} + \beta(v))^2$.
Then $x \maps \rho(x,v)$ $g_v$-concave on $\mathbb{R}^n$.

\label{lem:rhogv}
\end{lemma}

\begin{proof}
We first note that $\rho(x,v) > 0$, since $\rho(x,v)$ is the
solution to $g(x,rLv) = 0$, and $z \maps g(x,z)$ is convex (as a
convex quadratic map in $z$) as well as $g(x,0) < 0$ so that by
continuity any solution (if any) to $g(x,rLv) = 0$ must satisfy $r
> 0$. Consequently the composition $g_v(\rho(x,v))$ is
well-defined. Moreover, $g_v$ is clearly (continuously)
differentiable on $(0,\infty)$. Now observe that
\begin{align*}
\frac{1}{\rho(x,v)} &= \frac{h(x,v)}{-\beta(v) + \sqrt{\beta(v)^2
- b h(x,v)}}
                    = \frac{\beta(v) + \sqrt{\beta(v)^2 - b h(x,v)}}{-b}.
\end{align*}
So that, $\frac{-b}{\rho(x,v)} - \beta(v) = \sqrt{\beta(v)^2 - b
h(x,v)} \geq 0$, i.e.,
\begin{equation*}
g_v(\rho(x,v)) = -(\frac{b}{\rho(x,v)} + \beta(v))^2 = -\beta(v)^2
+ bh(x,v).
\end{equation*}
Now, $b < 0$ and $x \maps h(x,v)$ is convex, therefore $x \maps
g_v(\rho(x,v))$ is concave. It remains to show that $g_v$ is
strictly increasing on an appropriate domain.
To this end, observe that
\begin{itemize}
\item if $\beta(v) \leq 0$, we have for $t > 0$, $g_v'(t) = 2
\frac{b}{t^2} (\frac{b}{t} + \beta(v))$. Now, since $b < 0$ and
$\beta(v) \leq 0$, it follows that $(\frac{b}{t} + \beta(v)) < 0$
and evidently $\frac{b}{t^2} < 0$ so we conclude that $g_v'(t) >
0$, as was to be shown. Note that $g_v$ maps to
$(-\infty,-\beta(v)^2)$ for $\beta(v) \leq 0$.
\item if $\beta(v) > 0$, we still have, for $t > 0$, $g_v'(t) = (2
\frac{b}{t^2} (\frac{b}{t} + \beta(v))$ and $\frac{b}{t^2} < 0$.
However $(\frac{b}{t} + \beta(v)) < 0$ if and only if $t <
\frac{b}{-\beta(v)}$, so that $g_v$ is strictly increasing on
$(0,\frac{b}{-\beta(v)}]$ in this case only. We will now establish
that $\rho(x,v) \in (0, \frac{b}{-\beta(v)})$ holds for such
$v$. To this end recall
\begin{align*}
\frac{1}{\rho(x,v)} &= \frac{\beta(v) + \sqrt{\beta(v)^2 - b
h(x,v)}}{-b} ~\Longleftrightarrow~ \rho(x,v) = \frac{-b}{\beta(v)
+ \sqrt{\beta(v)^2 - b h(x,v)}},
\end{align*}
which with $\beta(v) + \sqrt{\beta(v)^2 - b h(x,v)} > \beta(v)
\geq 0$, gives $\rho(x,v) < \frac{b}{-\beta(v)}$ as desired.
\end{itemize}\hfill
\end{proof}

We can identify explicitly the inverse function of $g_v$, but for
this we will need to make a case-distinction.
\begin{itemize}
\item For $v$ such that $\beta(v) \leq 0$, $\reciprocal{g_v}
\fc{(-\infty,-\beta(v)^2)}{(0,\infty)}$ is given by
$\reciprocal{g_v}(t) = \frac{-b}{\sqrt{-t} + \beta(v)}$
\item For $v$ such that $\beta(v) > 0$, $\reciprocal{g_v}
\fc{(-\infty,0]}{(0, \frac{b}{-\beta(v)}]}$, is given by
$\reciprocal{g_v}(t) = \frac{-b}{\sqrt{-t} + \beta(v)}$.
\end{itemize}

We now employ Lemma \ref{lem:G1impliesG2} in order to combine
specialized concavity of $\rho$ with more usual notions.

\begin{lemma}
Let $v \in \sph$ be given and consider the map
$\reciprocal{g_v}(t) = \frac{-b}{\sqrt{-t} + \beta(v)}$ defined
from $(-\infty,-\beta(v)^2)$ to $(0,\infty)$ when $\beta(v) \leq
0$ and from $(-\infty,0]$ to $(0, \frac{b}{-\beta(v)}]$ otherwise.
This map is $-3$-concave on the set $(-\infty, -\beta(v)^2]$.
\label{lem:gv-3}
\end{lemma}

\begin{proof}
We will show this by a direct computation. In order to do so, let
$\alpha > -1$ be given but arbitrary and let $\psi$ denote the map
$t \mapsto (\frac{1}{\reciprocal{g_v}(t)})^{1+\alpha}$. The map
$\psi$ is clearly twice differentiable on the appropriate domain
and hence establishing the requested generalized concavity of
$\reciprocal{g_v}$ amounts to establishing convexity of $\psi$.
Now we establish
\begin{align*}
\psi(t) &=
(\frac{-1}{b})^{1+\alpha}(\sqrt{-t}+\beta(v))^{1+\alpha} \\
\psi'(t) &= -(\frac{-1}{b})^{1+\alpha} \frac{1+\alpha}{2}
(-t)^{\frac{-1}{2}} (\sqrt{-t} + \beta(v))^{\alpha} \\
\psi''(t) &= 
          (\frac{1+\alpha}{4
          b^2})(-t)^{-\frac{3}{2}}(\frac{1}{\reciprocal{g_v}(t)})^{\alpha - 1}\left (
          (\alpha-1)\sqrt{-t} - \beta(v) \right ).
\end{align*}
It is now clear that for $\beta(v) \leq 0$ and $\alpha > 1$, it
holds $\psi''(t) > 0$, implying the convexity of $\psi$, i.e., the
$-1-\alpha$-concavity of $\reciprocal{g_v}$. When $\beta(v) > 0$
holds, then for $\alpha \geq 2$, we observe that
$(\alpha-1)\sqrt{-t} - \beta(v) \geq \sqrt{-t} - \beta(v)$. The
latter is evidently positive whenever $t \leq -\beta^2(v)$. \hfill
\end{proof}

\begin{remark}
Let us observe that for any $v$ and $x$ with $g(x,0) < 0$,
$g_v(\rho(x,v)) \in (-\infty, -\beta(v)^2]$ holds. Indeed,
let us set $-t := -g_v(\rho(x,v)) = (\frac{b}{\rho(x,v)} +
\beta(v))^2$. Then as observed already: $\sqrt{-t} =
-\frac{b}{\rho(x,v)} - \beta(v) = \sqrt{\beta(v)^2 - bh(x,v)} \geq
0$. And consequently
\begin{equation*}
\sqrt{-t} - \beta(v) = -\frac{b}{\rho(x,v)} - 2\beta(v) =
\frac{1}{\rho(x,v)}( -b - 2\beta\rho(x,v) ).
\end{equation*}
But, $\rho(x,v)$ is the unique (positive) solution (in $r$) to the
equation $h(x,v) r^2 + 2\beta(v) r^2 + b = 0$. Therefore we have
$\sqrt{-t} - \beta(v) = \frac{1}{\rho(x,v)}h(x,v)\rho(x,v)^2 =
h(x,v)\rho(x,v) \geq 0$, thus establishing the claim.
\label{rem:rhoingv}
\end{remark}

We can now provide an eventual convexity statement for a
probability function involving$g$ defined in
\eqref{eq:gquad}.

\begin{proposition}
Consider the probability function $\phi(x) := \prb[g(x,\xi) \leq
0]$, where $x \in \Re^n$ is given and $g$ defined as in
\eqref{eq:gquad}. Let $\xi$ taking values in $\Re^m$ be an
elliptically symmetrically distributed random vector with mean $0$
and covariance matrix $\Sigma$ and associated radial distribution
$F_{\scr{R}}$. Let $F_{\scr{R}}$ be concave-$(-3)$ on the set $(0,
(t^\ast)^{-3}]$, with $t^*$ given by the first assumption of
Theorem\;\ref{thm:evencvx}.

Then for any $q \in (0,\frac 12)$, the set $M(p)$ is convex
provided that $p \geq p^\ast$ with
\begin{equation*}
p^\ast := \left(\frac 12 -
q\right)F_{\scr{R}}\left(\frac{t^\ast}{\delta(q)}\right) + \frac
12 + q,
\end{equation*}
where $\delta(q)$ is as in Theorem \ref{thm:evencvx}.

\label{prop:quadg}
\end{proposition}

\begin{proof}
The set $C := \col{x \in \Re^n \;: \rho(x,v) \geq t^\ast \foral v
\in \sph}$ is a convex set since $\rho(x,v)$ is quasi-concave in
$x$ as a result of convexity of $g$ in $x$. Moreover evidently
requisite 1. of Theorem \ref{thm:evencvx} holds. Next, let us
turn our attention to establishing requisite 2.

To this end, let $v \in \sph$ be given and let $g_v
\fc{\Re_+}{\Re_-}$ be the map defined in Lemma \ref{lem:rhogv}.
From this Lemma, we know in particular that $\rho$ is
$g_v$-concave on the range of $\rho(.,v)$. Hence in particular on
$C$. Let us thus pick $x_1,x_2 \in C$, $\lambda \in [0,1]$
arbitrarily. We thus obtain the estimate:
\begin{equation*}
\rho(\lambda x_1 + (1-\lambda)x_2,v) \geq \reciprocal{g_v}(
\lambda g_v(\rho(x_1,v)) + (1-\lambda)g_v(\rho(x_2,v)) ).
\end{equation*}
Moreover, by Remark \ref{rem:rhoingv} $g_v(\rho(x_1,v)) \in
(-\infty, -\beta(v)^2]$ (and likewise for $x_2$) and therefore
similarly $\lambda g_v(\rho(x_1,v)) + (1-\lambda)g_v(\rho(x_2,v))
\in (-\infty, -\beta(v)^2]$. Now, we may apply Lemma
\ref{lem:gv-3} to derive the further estimate:
\begin{eqnarray*}
\reciprocal{g_v}( \lambda g_v(\rho(x_1,v)) +
(1-\lambda)g_v(\rho(x_2,v)) ) &\geq&
m_{-3}(\reciprocal{g_v}(g_v(\rho(x_1,v))),\reciprocal{g_v}(g_v(\rho(x_2,v))),\lambda)\\
&=& m_{-3}(\rho(x_1,v), \rho(x_2,v), \lambda).
\end{eqnarray*}
Finally, since $F_{\scr{R}}$ is increasing as a distribution
function, we may pursue our estimates to obtain:
\begin{equation*}
F_{\scr{R}}(\rho(\lambda x_1 + (1-\lambda)x_2,v)) \geq
F_{\scr{R}}( m_{-3}(\rho(x_1,v), \rho(x_2,v), \lambda) ) \geq
\lambda F_{\scr{R}}( \rho(x_1,v) ) + (1-\lambda)F_{\scr{R}},
\end{equation*}
since $z \maps F_{\scr{R}}(z^{-\frac 13})$ is concave on the set
$(0, (t^\ast)^{-3}]$ by assumption and $\rho(x_1,v) \geq t^\ast$
(likewise for $x_2$) so that $\lambda \rho(x_1,v)^{-3} +
(1-\lambda) \rho(x_2,v)^{-3} \in (0, (t^\ast)^{-3}]$. This
concludes the proof of requisite 2.

Moreover note furthermore that the map $g$ defined in
\eqref{eq:gquad} is convex in the second argument $z$.
Consequently, for all $x \in M(\frac 12)$, we have $g(x,0) < 0$
and requisite 3. of Theorem \ref{thm:evencvx} holds for $p_0 =
\frac 12$. The resulting convexity now follows from applying
Theorem \ref{thm:evencvx}.
\end{proof}

\begin{remark}[Refined threshold]\label{rem:improve}
When $\xi$ is multi-variate Gaussian, the value $t^*$ can be
established to be $\sqrt{m+3}$ and moreover the threshold can be
strengthened to $p^\ast = \Phi(\sqrt{m+3})$. Note moreover that,
in all cases, we do not only have convexity of $M(p)$ for $p >
p^\ast$, but even concavity of $\phi$ on this set, which is a
rather strong property. Concavity of $\phi$ is not required to
derive convexity of its level-sets. Indeed, quasi-concavity would
suffice.

In view of this, it is of interest to recall results for the
specifically structured mapping $g(x,z) = z - x$. In this case,
$\rho$ can be proved to be concave \cite[Example
3.3]{vanAckooij_Malick_2017}). Using Theorem \ref{thm:evencvx}, we
would derive a stronger property: the concavity of $\phi$ on the sets of
the form $\col{x \in \Re^m \; : x \geq \mu + \sqrt{m-1}\norm{L}e}$
(with $e$ the all-one vector); if $\xi$ is Gaussian, this is
already known (e.g., \cite[Theorem 2.1]{Prekopa_2001}). Yet, it is
also well-known that multivariate Gaussian distribution functions
are log-concave (\cite{lagoa2005probabilistically}) and thus have all level sets convex. As a result,
the true threshold in that situation would be $p^*=0$.\hfill\qedsymbol
\end{remark}

\begin{figure}[!htbp]
\centering
\includegraphics[width=2.5in]{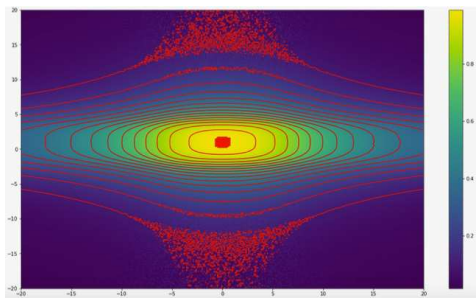}
\caption[caption]{Plot of the probability function $\phi$ for the
given quadratic problem.\label{fig:quad}}
\end{figure}

\begin{example}
As a numerical illustration, we propose to take both the decision
and random vector in dimension $2$, i.e., $m = n = 2$. We let
$\xi$ be a Gaussian random vector with mean $\mu = 0$ and
covariance matrix  $\Sigma = \begin{pmatrix}
0.01125 & 0.00675 \\
0.00675 & 0.2025
\end{pmatrix}$.
The mapping $W : \mathbb{R_+}^2  \rightarrow S_2^{++}$ is given by
$W(x_1,x_2)=
\begin{pmatrix}
x_1^2 + 0.5 & 0 \\
0 & |x_2-1|^3 + 0.2
\end{pmatrix}$, $2(a_1w_1^\mtr + a_2w_2^\mtr) = (1,1)$ and $b = -1$.

We obtain for these inputs the contour plot of
Figure\;\ref{fig:quad}. The contour lines are regularly drawn from
the probability value $0.1$ to $0.95$. We see on this graph how
the level sets tend to get convex as the level increases. The red
region in the center is the region where $\rho(x,v) \geq
\sqrt{m+3} = \sqrt{5}$ for any $v \in \mathbb{S}^{m-1}$.
Proposition \ref{prop:quadg} allows us to establish convexity for
$p^\ast = \Phi(\sqrt{5}) = 0.9873$. We see here that the set
obtained, is far, from being the largest one where convexity seems
to be present. This illustrates the discrepancy between the
current available threshold $p^*$, depending, in this case,
adversely on dimension and the practical situation, which seems to
exhibit convexity for thresholds significantly lower!\hfill\qedsymbol
\end{example}

We gather in a Python toolbox many useful functions for numerical
illustrations in this framework. In particular, we provide tools
for computing the probability function $\phi$ (with associated
function $\rho$) and plotting of its level-sets; as in Figure~\ref
{fig:quad}. We call this toolbox {\texttt pychance} and make it
publicly available on GitHub at
\url{https://github.com/yassine-laguel/pychance}.

\subsection{Examples with copul\ae\/}
\label{sec:excop}

We provide here eventual convexity statements based on results of Section \ref{subsec:copulae} for examples that are not covered by prior results. We choose these examples with decision vectors of dimension 2 to keep calculus simple without misrepresenting our results; this dimension does not impact our results, valid in infinite dimension. 

\begin{example}[Convexity statement for a $m$-dimensional copula]\label{ex:finfconcave}
Consider any  $m$ dimensional Archimedean copula, and construct
the probability functions
\[
\phi(x)= \C(F_1(h_1(x)),\ldots,F_m(h_m(x))),
\]
with a $G$-concave $h_1$ and $\alpha$-concave $h_i$ ($i\geq2$)
defined as follows. Given remark \ref{ex:relation}, we know that
$\C$ is $\invcon{G_1}$. Let $h_1 : \Re^2_{++} \rightarrow
\Re_{++}$ be defined by $h_1(x,y) = (\frac{y}{x^2})^2$ which is
neither concave nor convex, but $G$-concave with $G(x) =
-1/\sqrt{x}$ as easily shown\footnote{We see that $G(x)$ is
strictly increasing. To verify that $(x,y) \mapsto G(h_1(x,y)) = -
\frac{x^2}{y}$ is concave, just compute its Hessian
\begin{equation*}
-\frac{2}{y^3}
\begin{pmatrix}
y^2 & -xy \\
-xy & x^2
\end{pmatrix}
\end{equation*}
which is is negative semi-definite by a direct computation
(negative trace and zero determinant).}.
%
%
Let $h_i : \Re^2_{++} \rightarrow \Re_{++}$ be defined from any
convex positive function $f_i$ by $h_i(x,y) = f_i(x,y)^{\al_i}$
with $\al_i\in [-6,0)$. Let us define now the marginal
distribution functions. We take $F_1(z) = 1 - \exp(-\lambda z)$
the exponential distribution function of parameter $\lambda > 0$;
we easily\footnote{Observe first that $\reciprocal{G}(z) = z^{-2}$
is convex and strictly increasing. Verify that the composition
$F_1(\reciprocal{G}(z))$ is concave on some subset when $z$ takes
arguments in $(-\infty,0)$, by a direct computation of its second
derivative which gives $2f_1(\reciprocal{G}(z))z^{-4}(3-2\lambda
z^{-2})$ with the density function $f_1(z)=\la\exp(-\la z)$. } see
that it is $\invconG$ on the interval $I =  [b, \infty)$ with
$b=\sqrt{2\lambda/3}$. We take $F_i$ to be the Rayleigh
distribution function with parameter $\sigma=1.5$ for all $i \geq
2$.

We can then apply Theorem\;\ref{thm:convexity} to get the
convexity of $M(p)$ for any $p>p^*$ with
\[
p^* = \max\col{F_1\!\left(\frac{3}{2\la}\right), \max_{i\geq
2}F_{i}\left(\sqrt{(1-\al_i/2)\sigma}\right)} = 0.7769
\]
since $F_1\!\!\left(\frac{3}{2\la}\right) = 0.7769$ and the max
over $i\geq 2$ does not exceed $0.7364$. Thus the obtained
threshold is fixed and independent of the dimension of the random
vector.\hfill\qedsymbol
%
%
%
%
\end{example}

\begin{remark}[Independence from the dimension]
The previous example gives a situation wherein the threshold $p^*$
does not depend on the dimension $m$ of the random vector. The
underlying reason is that all mappings $h_i$ are generalized
concave with a set of parameters that do not degenerate with $m$.
It should also be observed that the involved parameters (see
e.g.,\;\cite[Table 1]{Henrion_Strugarek_2008}) depend continuously
on the parameters of the chosen marginal distributions and the
generalized concavity parameter of $h$. Therefore, should these
parameters belong to a compact set, the threshold $p^*$ does not
depend on dimension either. It is reasonable to assume that the
mappings $h_i$, marginal distributions, and their parameters are
homogenous. Typically the index $i$ expresses time and the
probability function stems from the desire to incorporate akin
features: each component $i$ is relatively similar.\hfill\qedsymbol
\end{remark}

\begin{example}[Convexity statement with functions from the literature]
Let us consider the functions of Example 1 of
\cite{Zadeh_Khorram_2012}: the two mappings from $\Re^2$ to
$\Re_+$
\[
h_1(x,y) = \exp{ -(x+y)^3 } \quad \text{and} \quad h_2(x,y) =
\frac{1}{x^2 + y^2 + 1}
\]
and the distributions $F_1 = \Phi$, the distribution function of a
standard normal Gaussian random variable and $F_2 = F_{\chi}$, the
distribution function of a $\chi$-random variable with $2$ degrees
of freedom. From Example\;\ref{ex:tbigalp}, the map $h_1$ can be
shown to be $G$-concave with a given strictly decreasing map $G$.
Moreover, Example \ref{ex:unusualG} gives us that $\Phi$ is
$\invconG$ on the set $(0, G(1.86)]$. Observe also that the map
$h_2$ is $(-1)$-concave (indeed, $\reciprocal{h_2}$ is convex) and
$F_{\chi}$ is $\invcon{(g_{-1})}$ (see Table\;1
in\;\cite{vanAckooij_Malick_2017}) on the set $(0,
\sqrt{3}^{-1}]$.

We extend now the situation of Example 1 of
\cite{Zadeh_Khorram_2012} which is restricted to the independent
copula. We consider here any $\invcon{g_1}$ copula (in particular
for all Archimedian copul\ae\/, by Example\;\ref{ex:relation}),
and we apply Theorem\;\ref{thm:convexity} to get that the feasible
set $M(p)$ with probability function
\begin{equation*}
\phi(x,y) = \C(F_1(h_1(x,y)), F_2(h_2(x,y))) ,
\end{equation*}
is convex for $p^\ast = \max\col{\Phi(1.86),
F_{\chi}(\sqrt{3})}=\max\col{0.9686, 0.7769}=0.9686$. Existing
theory cannot be applied to this case: \cite{Zadeh_Khorram_2012}
requires independent copula, and \cite{vanAckooij_2013} requires
$\alpha$-concavity (but $h_1$ is not $\alpha$-concave for any
$\alpha$).\hfill\qedsymbol
\end{example}

\begin{example}[Improved previous result for specific copula]
Let us go further with the previous example when $\C$ is
$\invcon{G_0}$ which subsumes the case of independent copula
considered in \cite{Zadeh_Khorram_2012}.
In this case, $z \maps \Phi(\reciprocal{G}(z))$ is log-concave on
the larger set $(0, G(t^\maltese)]$, with $t^\maltese = 1.6422$
(this value is numerically identified by employing the principles
of \cite[section 5.2]{vanAckooij_2013}), and $z \maps
F_{\chi}(z^{-1})$ is log-concave on an interval of the form
$(0,0.7563]$ (by Corollary\;\ref{cor:hs}). Then for any
$\invcon{G_0}$ copula (e.g., independent, Gumbel, Clayton),
Theorem\;\ref{thm:convexity} gives the threshold
\[
p^* = \max\col{\Phi(t^\maltese), F_{\chi}(1.3223)} =
\max\col{\Phi(1.6422), 0.5828} = 0.9497.
\]
This value is better than the threshold 
$\Phi(3)=0.9987$ given in \cite[Example 1]{Zadeh_Khorram_2012} for
the independent case only.\hfill\qedsymbol
\end{example}

\section*{Concluding Remarks}

In this paper we have provided general conditions under which
probabilistic constraints define a convex set. We have
investigated two different structures of the probability
functions: (i) a non-separable structure with elliptically random
vectors; (ii) a separable structure with the dependency expressed
by a given copula. In these situations, we have employed the
notion of $G$-concavity to reveal the level of underlying
convexity in the functions defining the constraints. The obtained
results are more general than those appearing in prior results.
The provided conditions can be verified from the nominal problem
data, as illustrated on various examples.
%
%

This work raises issues about the application of our convexity
results in practice for general chance-constrained optimization
problems (non-emptiness and convexity of feasible sets, guarantees
of solutions... see the discussion in Section \ref{sec:methods}).
A related question is the improvement of the computed thresholds.
In some specific situations indeed, it might be possible to lower
the computed threshold $p^*$ by exploiting structure of the
problems. More generally, more work is required in order to reduce
to gap between observed convexity and guaranteed convexity.

\section*{Acknowledgments}
We would like to acknowledge the partial financial support of PGMO
(Gaspard Monge Program for Optimization and operations research)
of the Hadamard Mathematic Foundation, through the project
``Advanced nonsmooth optimization methods for stochastic
programming".

\bibliographystyle{siam}
\bibliography{optim}


\end{document}